\documentclass[10pt,reqno,a4paper]{amsart}

\usepackage[T1]{fontenc}
\usepackage[utf8]{inputenc}
\usepackage{lmodern}
\usepackage{amsmath,amssymb,mathtools,mathrsfs}
\usepackage{bm}
\usepackage{enumitem}
\usepackage{microtype}
\usepackage[a4paper,top=24mm,bottom=25mm,left=25mm,right=25mm,headheight=13pt,headsep=7mm,footskip=10mm]{geometry}
\usepackage{graphicx}
\usepackage{float}
\usepackage{needspace}
\usepackage{placeins}
\usepackage[colorlinks=true,linkcolor=blue,citecolor=blue,urlcolor=blue,hypertexnames=false]{hyperref}

\hypersetup{
  pdftitle={Sharp History-Length Dependence of Mesh-Uniform Exponential Decay Rates for a Boundary-Memory Euler-Bernoulli Beam},
  pdfauthor={Gossrin Jean-Marc Bomisso},
  pdfsubject={Sharp history-length dependence and mesh-uniform stability for a memory-damped Euler-Bernoulli beam},
  pdfkeywords={optimal exponential decay, memory truncation, Euler-Bernoulli beam, Hermite-P1 approximation}
}

\allowdisplaybreaks
\numberwithin{equation}{section}

\newtheorem{theorem}{Theorem}[section]
\newtheorem{proposition}[theorem]{Proposition}
\newtheorem{lemma}[theorem]{Lemma}
\newtheorem{corollary}[theorem]{Corollary}

\theoremstyle{definition}

\theoremstyle{remark}
\newtheorem{remark}[theorem]{Remark}

\newcommand{\R}{\mathbb R}
\newcommand{\C}{\mathbb C}

\newcommand{\Xk}{\mathcal X_{\kappa}}
\newcommand{\Ak}{\mathbb A_{\kappa}}

\newcommand{\Ran}{\operatorname{Ran}}
\newcommand{\Real}{\operatorname{Re}}

\newcommand{\norm}[1]{\left\lVert #1\right\rVert}

\newcommand{\ip}[2]{\left\langle #1,#2\right\rangle}
\newcommand{\M}{\mathcal M}

\newcommand{\Ge}{\mathcal G}

\newcommand{\Xks}{\mathcal X_{\kappa,S}}
\newcommand{\Aks}{\mathbb A_{\kappa,S}}
\newcommand{\Ms}{\mathcal M_S}

\title[History-length dependence of mesh-uniform exponential decay]{Sharp History-Length Dependence of Mesh-Uniform Exponential Decay Rates for a Boundary-Memory Euler--Bernoulli Beam}

\author[G.~J.-M.~Bomisso]{Gossrin Jean-Marc Bomisso}
\makeatletter
\g@addto@macro\@setauthors{\par\vspace{1mm}{\centering\normalfont\small Laboratoire de Math\'ematiques et Informatique, UFR Sciences Fondamentales et Appliqu\'ees\\ Universit\'e Nangui Abrogoua, Abidjan, C\^ote d'Ivoire\\ \texttt{bogojm@yahoo.fr}\par}}
\makeatother

\makeatletter
\@ifundefined{subjclassname@2020}{\@namedef{subjclassname@2020}{2020 Mathematics Subject Classification}}{}
\makeatother

\begin{document}

\begin{abstract}
For a boundary-memory Euler--Bernoulli beam with polynomial relaxation, we quantify the loss of exponential stability as the retained history length $S$ grows. The continuous and conforming Hermite--$\mathbb P_1$ generators have zero-frequency resolvent norms bounded above and below by constant multiples of $S$. The bounds are uniform in the beam mesh and in every history mesh of step at most a fixed $h_0$. Their semigroup norms retain a plateau of duration proportional to $S$. Under positive distributed and boundary damping, matching exponential estimates hold with the optimal time scale $1+S$, uniformly in the meshes. Optimality refers to the semigroup norm on the augmented state space, which includes the history. Without distributed damping, the same time scale holds for the continuous truncated generator under a boundary condition on $\delta^{(3)}$. We also prove strong convergence of the finite-history semigroups to the full-history semigroup, uniformly on compact time intervals, under restriction and zero extension of histories. Remote-history packets give residuals of order $S^{-1}$ for the truncated generator and of order $S^{-1}+h$ after discretization, improved to $S^{-1}$ at zero frequency. Numerical tests confirm the residual orders.
\end{abstract}

\maketitle

\section{Introduction}
\label{sec:introduction}

The past-history formulation of Dafermos \cite{Dafermos1970} represents memory equations as semigroups on an enlarged state space, in which the history variable is transported by a weighted right translation. The spectral theory of this translation and its role in viscoelasticity are well established \cite{ChepyzhovPata2006,Pata2009,ContiDellOroPata2022,DellOroLiPata2026}. Dell'Oro, Li and Pata \cite{DellOroLiPata2026} analyze the translation generator itself, and Conti, Dell'Oro and Pata \cite{ContiDellOroPata2022} compare the history translation with the full semigroup. When the memory kernel decays slowly, histories concentrated far from the present are almost undamped. Such remote-history packets obstruct uniform decay of the coupled system. They become a quantitative issue as soon as the history is truncated to a finite interval $(0,S)$, which is unavoidable in computations.

Approximation of infinite delay by finite delay was considered by Conti, Marchini and Pata \cite{ContiMarchiniPata2012}. They recover infinite delay as a limit of finite delay on a large time scale, with a quantitative control of the discrepancy. Antonietti, Liverani and Pata \cite{AntoniettiLiveraniPata2023} studied numerical manifestations of non-superstable trajectories in viscoelasticity. Fabiano and Ito \cite{FabianoIto1990} gave a semigroup formulation and a convergence analysis for finite-history approximation schemes. Liu and Zheng \cite{LiuZheng1998} proved uniform exponential stability for the approximation schemes of Banks and Burns, Ito and Kappel, and Fabiano and Ito, under conditions on the kernel. Uniform stabilization of approximations of damped second-order systems was studied in \cite{ErvedozaZuazua2009,AbdallahNicaiseValeinWehbe2013}, for different operators or fixed memory regimes. The question addressed here is the sharp deterioration of the decay rate as $S\to\infty$, with bounds independent of the spatial and memory meshes.

We consider an Euler--Bernoulli beam clamped at one end, with boundary memory at the other end and kernel $\kappa(t)=b(1+t)^{-p}$, with $b>0$ and $p>1$. The history is truncated to $(0,S)$ and the system is approximated by cubic Hermite elements for the beam and continuous $\mathbb P_1$ elements for the history. Corrected remote packets give residuals of order $S^{-1}$ for the truncated generator, of order $S^{-1}+h$ after discretization, and of order $S^{-1}$ at zero frequency uniformly in both meshes. We deduce the sharp resolvent growth
\[
 cS\le\|\Aks^{-1}\|\le CS,
\]
with the same bounds for the Galerkin generator, and a plateau of the semigroup norms lasting a time proportional to $S$. Under strictly positive distributed and boundary damping, the matching estimate $Me^{-ct/(1+S)}$ holds for the truncated and Galerkin systems, uniformly in the meshes. A boundary multiplier gives the same optimal time scale for the continuous truncated system without distributed damping, and a separate result covers pure tip-velocity feedback. We also prove strong convergence of the truncated semigroups to the complete-history semigroup on bounded time intervals.

Decay of compatible physical trajectories with zero prehistory is treated in \cite{AnasseBomissoKpata2026}. The Hermite analysis of \cite{KoffiAnasseDansouBomisso2026} established a mesh-independent exponential estimate for a damped beam without truncated memory. The obstruction studied here comes from remote histories and not from high-frequency beam modes, which create other difficulties in numerical stabilization \cite{LeonZuazua2002,ErvedozaZuazua2009}. Section~\ref{sec:model} introduces the beam and its complete-history realization. Section~\ref{sec:truncation} treats truncation, discretization, resolvent growth and the optimal time scale. Appendix~\ref{app:static-transfer} computes the static history-to-beam transfer, with a critical exponent $p=2$, and Appendix~\ref{app:high-frequency} records a conditional high-frequency obstruction.

\section{The boundary-memory Euler--Bernoulli beam}
\label{sec:model}

\subsection{The model and the history space}

Throughout the paper, the relaxation kernel is given by
\begin{equation}
 \kappa(t)=b(1+t)^{-p},\quad\text{with}\quad b>0\ \text{and}\ p>1,
 \label{eq:kernel}
\end{equation}
and the memory density is defined by
\begin{equation*}
 m(s)=-\kappa'(s)=pb(1+s)^{-p-1}.
\end{equation*}

Let $L>0$. We consider the Euler--Bernoulli beam clamped at $x=0$, with boundary memory feedback at $x=L$, given by the system
\begin{equation}
 \begin{cases}
 \dfrac{\partial^2w}{\partial t^2}+\dfrac{\partial^4w}{\partial x^4}+\zeta\dfrac{\partial w}{\partial t}=0, & x\in(0,L),\ t>0,\\[2mm]
 w(0,t)=\dfrac{\partial w}{\partial x}(0,t)=0, & t>0,\\[2mm]
 \dfrac{\partial^2w}{\partial x^2}(L,t)+\delta^{(1)}\dfrac{\partial w}{\partial t}(L,t)
 +\delta^{(2)}\dfrac{\partial^2w}{\partial t\,\partial x}(L,t)
 +\delta^{(3)}\dfrac{\partial w}{\partial x}(L,t)=0, & t>0,\\[2mm]
 \dfrac{\partial^3w}{\partial x^3}(L,t)
 -\gamma\Bigl(\dfrac{\partial w}{\partial t}(L,t)+\displaystyle\int_0^\infty m(s)z(t,s)\,ds\Bigr)=0, & t>0,\\[2mm]
 \dfrac{\partial z}{\partial t}(t,s)+\dfrac{\partial z}{\partial s}(t,s)=\dfrac{\partial w}{\partial t}(L,t), & s>0,\ t>0,\\[2mm]
 z(t,0)=0, & t>0.
 \end{cases}
 \label{eq:beam-system}
\end{equation}
Here $\gamma>0$ is the boundary-feedback gain, denoted $A_b$ in \cite{AnasseBomissoKpata2026}. The coefficients $\zeta$, $\delta^{(1)}$, $\delta^{(2)}$ and $\delta^{(3)}$ are real.
We impose the following structural conditions:
\begin{equation}
 \zeta\ge0,\qquad \delta^{(1)}\ge0,\qquad \delta^{(3)}\ge0,\qquad \gamma>0
 \quad\text{and}\quad
 \delta^{(2)}\ge \frac{(\delta^{(1)})^2}{4\gamma}.
 \label{eq:structural}
\end{equation}
These assumptions guarantee generation and the estimates below. The Hermitian matrix
\begin{equation*}
 Q_L=
 \begin{pmatrix}
 \gamma & \delta^{(1)}/2\\
 \delta^{(1)}/2 & \delta^{(2)}
 \end{pmatrix}
\end{equation*}
is positive semidefinite, since
\[
 \det Q_L
 =\gamma\delta^{(2)}-\frac{(\delta^{(1)})^2}{4}\ge0.
\]
It is positive definite whenever the last inequality in
\eqref{eq:structural} is strict. In quadratic expressions,
$\mathbf v^*$ denotes the conjugate transpose of a complex vector
$\mathbf v$.

We write $v=w_t$ for the beam velocity. We define the spaces
\[
 H_E^2(0,L)=\{w\in H^2(0,L):w(0)=w'(0)=0\}
\]
and
\[
 L_m^2(\R_+)=
 \left\{z:\R_+\to\C:\int_0^\infty m(s)|z(s)|^2\,ds<\infty\right\}.
\]
Since
\begin{equation*}
 \int_0^\infty m(s)\,ds=b,
\end{equation*}
the functional
\begin{equation*}
 Jz:=\int_0^\infty m(s)z(s)\,ds
\end{equation*}
is continuous on $L_m^2(\R_+)$ and satisfies
\begin{equation}
 |Jz|\le \sqrt b\,\norm{z}_{L_m^2}.
 \label{eq:Jbound}
\end{equation}
We denote this weighted history space by $\M=L_m^2(\R_+)$. Its right-translation generator is
\[
 \Ge z=-z_s,\quad\text{with}\quad 
 D(\Ge)=\{z\in\M: z\text{ is locally absolutely continuous},\ z_s\in\M,\ z(0)=0\}.
\]
The translation is contractive because $m$ is nonincreasing.

The augmented energy space is
\begin{equation*}
 \Xk=H_E^2(0,L)\times L^2(0,L)\times L_m^2(\R_+),
\end{equation*}
with inner product, linear in the first argument,
\begin{align}
 \ip{(w,v,z)}{(\phi,\psi,\eta)}_{\Xk}
 ={}&\int_0^L w''\overline{\phi''}\,dx
 +\delta^{(3)}w'(L)\overline{\phi'(L)}
 +\int_0^L v\overline\psi\,dx \notag\\
 &+\gamma\int_0^\infty m(s)z(s)\overline{\eta(s)}\,ds.
 \label{eq:inner-product}
\end{align}
The norm is equivalent to the natural product norm because $\norm{w''}_{L^2}$ is a norm on $H_E^2(0,L)$.

System \eqref{eq:beam-system} is written as $U_t=\Ak U$ for $U=(w,v,z)$, where the generator is defined by
\begin{equation*}
 \Ak(w,v,z)=\bigl(v,-w^{(4)}-\zeta v,\,v(L)-z_s\bigr),
\end{equation*}
with domain
\begin{equation}
 \begin{aligned}
 D(\Ak)=\bigl\{(w,v,z)\in\Xk:{}
 &w\in H^4(0,L)\cap H_E^2(0,L),\\
 &v\in H_E^2(0,L)\quad\text{and}\quad z\in D(\Ge),\\
 &w''(L)+\delta^{(1)}v(L)+\delta^{(2)}v'(L)+\delta^{(3)}w'(L)=0,\\
 &w'''(L)-\gamma\bigl(v(L)+Jz\bigr)=0\bigr\}.
 \end{aligned}
 \label{eq:domain}
\end{equation}

\subsection{Dissipativity and generation}

\begin{proposition}\label{prop:dissipativity}
Under \eqref{eq:structural}, the operator $\Ak$ is dissipative on $\Xk$.
\end{proposition}

\begin{proof}
Let $U=(w,v,z)\in D(\Ak)$. Integrating twice by parts in the beam component and using the two boundary conditions gives
\begin{align*}
 \Real\ip{\Ak U}{U}_{\Xk}
 ={}&-\zeta\norm{v}_{L^2}^2
 -\gamma|v(L)|^2
 -\delta^{(1)}\Real\bigl(v(L)\overline{v'(L)}\bigr)
 -\delta^{(2)}|v'(L)|^2\\
 &-\gamma\Real\bigl(Jz\,\overline{v(L)}\bigr)
 +\gamma\Real\int_0^\infty m(s)v(L)\overline{z(s)}\,ds\\
 &-\gamma\Real\int_0^\infty m(s)z_s(s)\overline{z(s)}\,ds.
\end{align*}
The two mixed memory terms cancel. Since
\[
 |m'(s)|=\frac{p+1}{1+s}m(s)\le(p+1)m(s),
\]
a cutoff argument justifies integration by parts without any pointwise condition at infinity. For cutoffs \(\chi_R\) equal to one on \([0,R]\) and supported in \([0,R+1]\), the cutoff boundary term is bounded by
\[
 C\int_R^{R+1}m(s)|z(s)|^2\,ds\to0,
\]
while, with $\mathbf 1_A$ the characteristic function of a set $A$, the remaining tail is controlled by
\[
 \|z_s\mathbf1_{(R,\infty)}\|_{L_m^2}
 \|z\mathbf1_{(R,\infty)}\|_{L_m^2}\to0.
\]
Hence
\[
 \Real\int_0^\infty m(s)z_s(s)\overline{z(s)}\,ds
 =-\frac12\int_0^\infty m'(s)|z(s)|^2\,ds.
\]
With
\[
 \mathbf v_L:=\begin{pmatrix}v(L)\\v'(L)\end{pmatrix},
\]
we obtain the Hermitian form
\begin{align*}
 \Real\ip{\Ak U}{U}_{\Xk}
 ={}&-\zeta\norm{v}_{L^2}^2-\mathbf v_L^{\,*}Q_L\mathbf v_L
 +\frac{\gamma}{2}\int_0^\infty m'(s)|z(s)|^2\,ds\le0.
\end{align*}
\end{proof}

\paragraph{The resolvent equation.}

Fix $\lambda>0$ and let $F=(f,g,\vartheta)\in\Xk$. We solve
\begin{equation}
 (\lambda I-\Ak)(w,v,z)=(f,g,\vartheta).
 \label{eq:resolvent}
\end{equation}
The first two equations give
\begin{equation}
 v=\lambda w-f
 \label{eq:v-elimination}
\end{equation}
and
\begin{equation}
 w^{(4)}+\lambda(\lambda+\zeta)w=F_\lambda,\quad\text{with}\quad F_\lambda:=g+(\lambda+\zeta)f.
 \label{eq:reduced-ode}
\end{equation}
The transport equation is
\begin{equation*}
 z_s+\lambda z=v(L)+\vartheta\quad\text{and}\quad z(0)=0.
\end{equation*}
We set
\begin{equation*}
 q_\lambda(s)=\frac{1-e^{-\lambda s}}{\lambda}
 \quad\text{and}\quad
 (R_\lambda \vartheta)(s)=\int_0^s e^{-\lambda(s-r)}\vartheta(r)\,dr.
\end{equation*}
The solution is given by
\begin{equation}
 z(s)=q_\lambda(s)v(L)+(R_\lambda \vartheta)(s).
 \label{eq:z-resolvent}
\end{equation}

\begin{lemma}\label{lem:transport-resolvent}
The operator $R_\lambda$ is bounded on $L_m^2(\R_+)$ and
\begin{equation*}
 \norm{R_\lambda \vartheta}_{L_m^2}\le \frac1\lambda\norm{\vartheta}_{L_m^2}.
\end{equation*}
Moreover, the function $z$ in \eqref{eq:z-resolvent} belongs to $D(\Ge)$.
\end{lemma}

\begin{proof}
Cauchy--Schwarz and Tonelli's theorem, together with the monotonicity of \(m\), give
\[
 \|R_\lambda \vartheta\|_{L_m^2}^2
 \le \frac1{\lambda}\int_0^\infty |\vartheta(r)|^2
 \left(\int_r^\infty m(s)e^{-\lambda(s-r)}\,ds\right)dr
 \le\frac1{\lambda^2}\|\vartheta\|_{L_m^2}^2.
\]
Moreover \(0\le q_\lambda\le\lambda^{-1}\), so \(z\in L_m^2\). Since \(m\) is bounded below on every compact interval, \(\vartheta\in L^2_{\mathrm{loc}}(\R_+)\). The explicit formula \eqref{eq:z-resolvent} therefore shows that \(z\) is locally absolutely continuous. From \(z_s=-\lambda z+v(L)+\vartheta\) and \(\int_0^\infty m=b\), one obtains \(z_s\in L_m^2\). Thus \(z\in D(\Ge)\).
\end{proof}

We define
\begin{equation*}
 \rho_\lambda:=Jq_\lambda>0\quad\text{and}\quad H_\lambda:=J(R_\lambda \vartheta).
\end{equation*}
Then
\begin{equation}
 Jz=\rho_\lambda v(L)+H_\lambda.
 \label{eq:Jz-resolvent}
\end{equation}
Substituting \eqref{eq:v-elimination} and \eqref{eq:Jz-resolvent} into the boundary conditions gives
\begin{align}
 w''(L)+\lambda\delta^{(1)}w(L)
 +(\lambda\delta^{(2)}+\delta^{(3)})w'(L)&=B_1,
 \label{eq:bc-reduced1}\\
 w'''(L)-\gamma\lambda(1+\rho_\lambda)w(L)&=B_2,
 \label{eq:bc-reduced2}
\end{align}
where
\begin{equation*}
 B_1=\delta^{(1)}f(L)+\delta^{(2)}f'(L)
 \quad\text{and}\quad
 B_2=\gamma\bigl(H_\lambda-(1+\rho_\lambda)f(L)\bigr).
\end{equation*}

For $w,\phi\in H_E^2(0,L)$, set
\begin{align*}
 a_\lambda(w,\phi)
 ={}&\int_0^L w''\overline{\phi''}\,dx
 +\lambda(\lambda+\zeta)\int_0^L w\overline\phi\,dx
 +\gamma\lambda(1+\rho_\lambda)w(L)\overline{\phi(L)}\notag\\
 &+\lambda\delta^{(1)}w(L)\overline{\phi'(L)}
 +(\lambda\delta^{(2)}+\delta^{(3)})w'(L)\overline{\phi'(L)},
\end{align*}
and
\begin{equation*}
 \Lambda_\lambda(\phi)
 =\int_0^L F_\lambda\overline\phi\,dx
 -B_2\overline{\phi(L)}+B_1\overline{\phi'(L)}.
\end{equation*}
Both are continuous on $H_E^2(0,L)$.

\begin{lemma}\label{lem:coercivity}
For every $\lambda>0$, the form $a_\lambda$ is coercive on $H_E^2(0,L)$.
\end{lemma}

\begin{proof}
Using Young's inequality in the form
\[
 \delta^{(1)}|ab|
 \le \gamma|a|^2+\frac{(\delta^{(1)})^2}{4\gamma}|b|^2,
\]
we obtain
\begin{align*}
 \Real a_\lambda(w,w)
 \ge{}&\norm{w''}_{L^2}^2
 +\lambda(\lambda+\zeta)\norm{w}_{L^2}^2
 +\lambda \gamma\rho_\lambda |w(L)|^2\notag\\
 &+\left[\lambda\left(\delta^{(2)}
 -\frac{(\delta^{(1)})^2}{4\gamma}\right)+\delta^{(3)}\right]|w'(L)|^2.
\end{align*}
By \eqref{eq:structural}, the last coefficient is nonnegative. Since $\norm{w''}_{L^2}$ is a norm equivalent to the $H^2(0,L)$ norm on $H_E^2(0,L)$, the form is coercive on $H_E^2(0,L)$.
\end{proof}

\begin{theorem}[Generation on the complete history space]\label{thm:generation}
Under \eqref{eq:structural}, the operator $\Ak$ is maximal dissipative on $\Xk$. Consequently, it generates a contraction $C_0$-semigroup $(T(t))_{t\ge0}$ on $\Xk$.
\end{theorem}

\begin{proof}
By Lemma~\ref{lem:coercivity} and the complex Lax--Milgram theorem, there is a unique $w\in H_E^2(0,L)$ satisfying
\begin{equation}
 a_\lambda(w,\phi)=\Lambda_\lambda(\phi),
 \quad\text{for all }\phi\in H_E^2(0,L).
 \label{eq:weak-resolvent}
\end{equation}
Equation \eqref{eq:reduced-ode} holds in distributions. Hence
$w^{(4)}\in L^2(0,L)$. The one-dimensional regularity implication
\[
 w\in H^2(0,L),\quad w^{(4)}\in L^2(0,L)
 \quad\Longrightarrow\quad w\in H^4(0,L)
\]
gives $w\in H^4(0,L)$. Integrating by parts in \eqref{eq:weak-resolvent} and using the independence of the traces $\phi(L)$ and $\phi'(L)$ recovers \eqref{eq:bc-reduced1}--\eqref{eq:bc-reduced2}.

Define $v$ by \eqref{eq:v-elimination} and $z$ by \eqref{eq:z-resolvent}. Then $v\in H_E^2(0,L)$ and, by Lemma~\ref{lem:transport-resolvent}, $z\in D(\Ge)$. The reduced boundary equations are exactly the original domain conditions after substitution. Thus $(w,v,z)\in D(\Ak)$ and solves \eqref{eq:resolvent}. Hence
\begin{equation}
 \Ran(\lambda I-\Ak)=\Xk.
 \label{eq:range}
\end{equation}

Since $\Xk$ is a Hilbert space, Proposition~\ref{prop:dissipativity}, \eqref{eq:range} and the Lumer--Phillips theorem \cite{Pazy} give the density of $D(\Ak)$ and the contraction semigroup.
\end{proof}

\section{History truncation and finite elements}
\label{sec:truncation}

\subsection{Finite-history truncation}
\label{subsec:finite-history-pseudo}

For a direct numerical treatment of the history variable, we truncate the interval to $(0,S)$. Fix $S>2$ and set
\[
 \Ms=L_m^2(0,S)
 \quad\text{and}\quad
 J_S z=\int_0^S m(s)z(s)\,ds.
\]
Let $\mathcal G_S z=-z_s$ on
\[
\begin{aligned}
 D(\mathcal G_S)=\{z\in\Ms:\;& z\text{ is locally absolutely continuous},\\
 &z_s\in\Ms\ \text{and}\ z(0)=0\}.
\end{aligned}
\]
We set
\[
 \Xks=H_E^2(0,L)\times L^2(0,L)\times\Ms
\]
with the energy inner product inherited from \eqref{eq:inner-product}, with the history integral restricted to $(0,S)$. The finite-history generator is
\[
 \Aks(w,v,z)=\bigl(v,-w^{(4)}-\zeta v,\,v(L)-z_s\bigr).
\]
Its domain is obtained from \eqref{eq:domain} by replacing $D(\Ge)$ by $D(\mathcal G_S)$ and $J$ by $J_S$. No boundary condition is imposed at the outflow point $s=S$.

\begin{proposition}[Finite-history generation]\label{prop:finite-history-generation}
Under \eqref{eq:structural}, $\Aks$ is maximal dissipative on $\Xks$ and hence generates a contraction $C_0$-semigroup.
\end{proposition}

\begin{proof}
For fixed \(S\), positivity and continuity of \(m\) on \([0,S]\) make the weighted norm equivalent to the usual \(L^2(0,S)\) norm. Hence
\[
 D(\mathcal G_S)=\{z\in H^1(0,S):z(0)=0\}
\]
is dense and \(\mathcal G_S=-\partial_s\) is closed. No condition is imposed at the outflow point \(s=S\). Integration by parts gives
\[
 -\Real\int_0^S m z_s\overline z\,ds
 =-\frac12m(S)|z(S)|^2+\frac12\int_0^S m'(s)|z(s)|^2\,ds\le0,
\]
so the dissipativity proof of Proposition~\ref{prop:dissipativity} carries over, with the additional nonpositive outflow term.

For \(\lambda>0\), the transport resolvent is
\[
 z(s)=q_\lambda(s)v(L)+(R_{\lambda,S}\vartheta)(s),\quad\text{with}\quad 
 (R_{\lambda,S}\vartheta)(s)=\int_0^s e^{-\lambda(s-r)}\vartheta(r)\,dr,
\]
and \(R_{\lambda,S}\) is bounded on \(\Ms\). Since
\[
 \rho_{\lambda,S}=\int_0^S m(s)q_\lambda(s)\,ds>0,
\]
the reduced beam problem is the same as in Lemma~\ref{lem:coercivity}, with \(\rho_\lambda\) replaced by \(\rho_{\lambda,S}\), and
\[
 \Real a_{\lambda,S}(w,w)\ge
 \|w''\|_{L^2}^2+
 \left[\lambda\left(\delta^{(2)}-\frac{(\delta^{(1)})^2}{4\gamma}\right)
 +\delta^{(3)}\right]|w'(L)|^2.
\]
The Lax--Milgram and regularity arguments of Theorem~\ref{thm:generation} show that
\(\lambda I-\Aks\) is onto \(\Xks\).
Thus $\Aks$ is maximal dissipative, and the Lumer--Phillips theorem gives the contraction semigroup.
\end{proof}

\begin{theorem}[Strong consistency of history truncation]
\label{thm:truncation-consistency}
For $S>2$, define $P_S:\Xk\to\Xks$ by retaining both mechanical components and restricting the history to $(0,S)$. Define $E_S:\Xks\to\Xk$ by zero extension of the history, leaving the mechanical components unchanged. Under \eqref{eq:structural}, for every $U\in\Xk$ and $T>0$,
\begin{equation}\label{eq:strong-truncation-consistency}
 \lim_{S\to\infty}\sup_{0\le t\le T}
 \|E_S e^{t\Aks}P_SU-e^{t\Ak}U\|_{\Xk}=0.
\end{equation}
The convergence holds for each fixed initial state, not in operator norm.
\end{theorem}
\begin{proof}
The restriction $P_S$ and zero extension $E_S$ are contractions and satisfy $P_SE_S=I_{\Xks}$. Since $m$ is integrable and $U\in\Xk$, we have $E_SP_SU\to U$ in $\Xk$.

Fix $\lambda>0$ and $F=(f,g,\vartheta)\in\Xk$. Write $U_\lambda=(\lambda I-\Ak)^{-1}F$ and $U_{\lambda,S}=(\lambda I-\Aks)^{-1}P_SF$. The transport resolvent formula \eqref{eq:z-resolvent} gives, with $q_\lambda(s)=(1-e^{-\lambda s})/\lambda$,
\[
 z_\lambda(s)=q_\lambda(s)v_\lambda(L)+R_\lambda\vartheta(s)\quad\text{and}\quad
 z_{\lambda,S}(s)=q_\lambda(s)v_{\lambda,S}(L)+R_\lambda\vartheta(s),\quad\text{for }0<s<S.
\]
Set $\rho_{\lambda,S}=\int_0^Smq_\lambda\,ds$ and $H_{\lambda,S}=\int_0^SmR_\lambda\vartheta\,ds$. Since $q_\lambda\in\M$ and $R_\lambda\vartheta\in\M$ by Lemma~\ref{lem:transport-resolvent}, Cauchy--Schwarz and the integrability of $m$ give
\[
 \rho_{\lambda,S}\longrightarrow\rho_\lambda\quad\text{and}\quad H_{\lambda,S}\longrightarrow H_\lambda.
\]
The reduced variational problems are \eqref{eq:weak-resolvent} for the full history and its truncated counterpart obtained by replacing $\rho_\lambda,H_\lambda$ with $\rho_{\lambda,S},H_{\lambda,S}$. Their coercivity constant is bounded below independently of $S$ by Lemma~\ref{lem:coercivity}, since all terms involving $\rho_{\lambda,S}$ are nonnegative. Subtraction of the two variational equations therefore gives $w_{\lambda,S}\to w_\lambda$ in $H_E^2(0,L)$. Consequently $v_{\lambda,S}=\lambda w_{\lambda,S}-f\to v_\lambda$ in $H_E^2(0,L)$, and
\[
 \|z_{\lambda,S}-P_Sz_\lambda\|_{\Ms}
 \le\|q_\lambda\|_{\M}\,|v_{\lambda,S}(L)-v_\lambda(L)|\longrightarrow0.
\]
Together with $E_SP_SU_\lambda\to U_\lambda$, this establishes the strong resolvent convergence
\[
 E_S(\lambda I-\Aks)^{-1}P_SF\longrightarrow
 (\lambda I-\Ak)^{-1}F \quad\text{in }\Xk.
\]
Both generators are maximal dissipative, by Theorem~\ref{thm:generation} and Proposition~\ref{prop:finite-history-generation}. To use Trotter--Kato on a common Hilbert space, set $Q_S=E_SP_S$, the orthogonal projection onto the zero-extended truncated space. Define
\[
 \widehat T_S(t)=E_S e^{t\Aks}P_S+(I-Q_S),\quad\text{for all }t\ge0.
\]
This is a contraction $C_0$-semigroup on $\Xk$ with identity at $t=0$. Its resolvent is
\[
 (\lambda I-\widehat A_S)^{-1}F
 =E_S(\lambda I-\Aks)^{-1}P_SF
 +\lambda^{-1}(I-Q_S)F.
\]
The first term converges to $(\lambda I-\Ak)^{-1}F$ and the second tends to zero, since $Q_SF\to F$ strongly. The contraction version of the Trotter--Kato theorem on $\Xk$ \cite{Pazy} yields
\[
 \sup_{0\le t\le T}\|\widehat T_S(t)U-e^{t\Ak}U\|_{\Xk}\longrightarrow0.
\]
Finally, $\widehat T_S(t)U-E_Se^{t\Aks}P_SU=(I-Q_S)U$, whose norm tends to zero. This proves \eqref{eq:strong-truncation-consistency}.
\end{proof}

For a closed densely defined operator $B$ and $\varepsilon>0$, we use the closed $\varepsilon$-pseudospectrum
\[
 \sigma_\varepsilon(B)
 =\sigma(B)\cup
 \bigl\{\lambda\in\rho(B):
 \|(\lambda I-B)^{-1}\|\ge\varepsilon^{-1}\bigr\},
\]
so that the existence of a unit vector $u\in D(B)$ with $\|(\lambda I-B)u\|\le\varepsilon$ forces $\lambda\in\sigma_\varepsilon(B)$. See, for example, \cite{TrefethenEmbree}.

\begin{lemma}[Cancellation of the remote-history moment]\label{lem:moment-cancellation}
Let $p>1$, $R_S=3S/4$, $\ell_S=S/8$ and $I_S=[5S/8,7S/8]$. Fix a nonzero
$\varphi\in C_c^\infty((-1,1))$. For each fixed $\beta\ne0$, let
\[
 \widetilde z_{S,\beta}(s)=a_{S,\beta}e^{-i\beta s}
 \varphi\!\left(\frac{s-R_S}{\ell_S}\right),\quad\text{with}\quad \|\widetilde z_{S,\beta}\|_{\Ms}=1.
\]
For each integer $k\ge0$ and each $S\ge4$,
\begin{equation}\label{eq:oscillatory-moment}
 |J_S\widetilde z_{S,\beta}|\le C_k S^{-p/2}
 \min\{1,(|\beta|\ell_S)^{-k}\},\quad\text{for all }\beta\ne0,
\end{equation}
where $C_k$ is independent of $S$ and $\beta$. In particular, for fixed
$\beta\ne0$ the moment is $O_{\beta,N}(S^{-N})$ for every $N>0$.

At $\beta=0$ there is a normalized function $z_{S,0}\in C_c^\infty(I_S)$
such that
\begin{equation}\label{eq:zero-moment-direct}
 J_Sz_{S,0}=0,\qquad
 \|\mathcal G_Sz_{S,0}\|_{\Ms}\le CS^{-1}\quad\text{and}\quad
 \|z_{S,0}^{(j)}\|_{\Ms}\le CS^{-j},\quad\text{for }j=0,1,2.
\end{equation}
The constants are independent of $S$.
\end{lemma}

\begin{proof}
For $\beta\ne0$, substitute $s=R_S+\ell_Sr$ and write
\[
 J_S\widetilde z_{S,\beta}
 =a_{S,\beta}\ell_Se^{-i\beta R_S}
 \int_{-1}^1g_S(r)e^{-i\beta\ell_Sr}\,dr,\quad\text{with}\quad g_S(r)=m(R_S+\ell_Sr)\varphi(r).
\]
For every integer $j\ge0$, the polynomial weight satisfies
$|m^{(j)}(s)|\le C_jm(s)(1+s)^{-j}$. Since
$\ell_S/(1+s)\le1/5$ on $I_S$, Leibniz' rule gives
$\|g_S^{(k)}\|_{L^1(-1,1)}\le C_km_S^+$, where
$m_S^\pm$ are the maximum and minimum of $m$ on $I_S$. The function
$\varphi$ vanishes near both endpoints. Thus repeated integration by parts
gives
\[
 \left|\int_{-1}^1g_S(r)e^{-i\beta\ell_Sr}\,dr\right|
 \le C_km_S^+\min\{1,(|\beta|\ell_S)^{-k}\}.
\]
Normalization and $m_S^+/m_S^-\le(7/5)^{p+1}$ imply
$a_{S,\beta}\ell_Sm_S^+\le C(m_S^-\ell_S)^{1/2}\le CS^{-p/2}$.
This proves \eqref{eq:oscillatory-moment}.

For $\beta=0$, select nonnegative nonzero functions
$\varphi_1,\varphi_2\in C_c^\infty((-1,1))$ with disjoint supports, and set
\[
 c_S=\frac{\int_{-1}^1m(R_S+\ell_Sr)\varphi_1(r)\,dr}
 {\int_{-1}^1m(R_S+\ell_Sr)\varphi_2(r)\,dr}\quad\text{and}\quad
 \Phi_S=\varphi_1-c_S\varphi_2.
\]
Uniform comparability of the weight on $I_S$ shows that
$0<c\le c_S\le C$, while the disjoint supports imply that the weighted
$L^2$ norm of $\Phi_S$ is bounded above and below, after factoring out
$m_S^-$. Choose $a_S>0$ so that
\[
 z_{S,0}(s)=a_S\Phi_S\!\left(\frac{s-R_S}{\ell_S}\right),\quad\text{with}\quad \|z_{S,0}\|_{\Ms}=1.
\]
By construction $J_Sz_{S,0}=0$. For $j=1,2$, a change of variables
and weight comparability give
$\|z_{S,0}^{(j)}\|_{\Ms}\le C\ell_S^{-j}$.
Since $\mathcal G_S=-\partial_s$ and $\ell_S=S/8$, this proves
\eqref{eq:zero-moment-direct}.
\end{proof}

\begin{proposition}[Pseudospectral persistence under history truncation]
\label{prop:finite-history-pseudo}
For each fixed $\beta\in\R$ there exist $S_0>0$ and $C_\beta>0$,
independent of $S$, such that for every $S\ge S_0$ there is
$U_{S,\beta}\in D(\Aks)$ with
\begin{equation}\label{eq:finite-history-residual}
 \|U_{S,\beta}\|_{\Xks}=1\quad\text{and}\quad
 \|(i\beta I-\Aks)U_{S,\beta}\|_{\Xks}\le C_\beta S^{-1}.
\end{equation}
Consequently $i\beta\in\sigma_{C_\beta S^{-1}}(\Aks)$ for all sufficiently
large $S$. For $\beta$ in a fixed compact subset of $\R\setminus\{0\}$,
the constant can be chosen uniformly in $\beta$.
\end{proposition}

\begin{proof}
If $\beta=0$, take $z_{S,\beta}=z_{S,0}$ from
Lemma~\ref{lem:moment-cancellation}. Its moment vanishes and its
residual is bounded by $CS^{-1}$.

Fix $\beta\ne0$. Take $\widetilde z_{S,\beta}$ from the same lemma.
Differentiation of the explicit packet gives
$\|(i\beta I-\mathcal G_S)\widetilde z_{S,\beta}\|_{\Ms}
=\ell_S^{-1}\|a_{S,\beta}e^{-i\beta s}\varphi'((s-R_S)/\ell_S)\|_{\Ms}\le CS^{-1}$,
by uniform comparability of $m$ on $I_S$.
Choose a nonnegative nonzero $\psi\in C_c^\infty((1,2))$, so that $J\psi>0$. Its support is disjoint from $I_S$
for large $S$. Set
\begin{equation}
 c_{S,\beta}=\frac{J_S\widetilde z_{S,\beta}}{J\psi},\qquad
 z_{S,\beta}^0=\widetilde z_{S,\beta}-c_{S,\beta}\psi\quad\text{and}\quad
 z_{S,\beta}=\frac{z_{S,\beta}^0}{\|z_{S,\beta}^0\|_{\Ms}}.
 \label{eq:corrected-packet}
\end{equation}
Then $J_Sz_{S,\beta}=0$. Disjoint supports give
$\|z_{S,\beta}^0\|_{\Ms}^2
 =1+|c_{S,\beta}|^2\|\psi\|_{\Ms}^2\ge1$.
Moreover, by Lemma~\ref{lem:moment-cancellation} with $k=1$,
$|c_{S,\beta}|\le C|\beta|^{-1}S^{-1-p/2}$.
Hence
$\|(i\beta I-\mathcal G_S)z_{S,\beta}\|_{\Ms}\le C_\beta S^{-1}$.
For either case set $U_{S,\beta}=\gamma^{-1/2}(0,0,z_{S,\beta})$.
The moment constraint implies $U_{S,\beta}\in D(\Aks)$ and
$\Aks U_{S,\beta}=\gamma^{-1/2}(0,0,\mathcal G_Sz_{S,\beta})$.
This proves \eqref{eq:finite-history-residual} and the pseudospectral
inclusion. The bounds are uniform on compact frequency sets avoiding zero.
\end{proof}

\begin{remark}[Finite-history translation]\label{rem:finite-transport-spectrum}
For every \(\lambda\in\C\), the equation \((\lambda I-\mathcal G_S)z=f\) has the Volterra solution
\[
 z(s)=\int_0^s e^{-\lambda(s-r)}f(r)\,dr,
\]
so \(\sigma(\mathcal G_S)=\varnothing\). Proposition~\ref{prop:finite-history-pseudo} gives the corresponding
pseudospectral estimate after truncation. For $\beta\ne0$, oscillation makes the moment correction negligible. At $\beta=0$, the two-profile construction in Lemma~\ref{lem:moment-cancellation} gives $O(S^{-1})$ for every $p>1$. A single packet with a fixed corrector gives only $O(S^{-1}+S^{-p/2})$. For fixed $S$, $\Aks$ has compact resolvent. Its resolvent maps into the operator domain, whose beam components lie in $H^4(0,L)$ and $H^2(0,L)$ and whose history component lies in $H^1(0,S)$. The corresponding embeddings into $H_E^2(0,L)$, $L^2(0,L)$ and $\Ms$ are compact. Thus the truncated operator has compact resolvent, while its
pseudospectral residuals approach zero as $S\to\infty$.
\end{remark}

\subsection{\texorpdfstring{Galerkin approximation with conforming Hermite and $\mathbb P_1$ spaces}{Galerkin approximation with conforming Hermite and P1 spaces}}

We use Galerkin spaces conforming in the energy space, but not in the domain of the continuous generator. The continuous $\mathbb P_1$ history space permits an exact energy identity after $L_m^2$ projection of the transport derivative, while the inflow condition is imposed strongly. This dissipativity property is essential to the mesh-uniform stability argument. For $k\ge0$, let
$\mathbb P_k$ denote the polynomials of degree at most $k$. Fix a beam mesh
\[
 0=x_0<x_1<\cdots<x_{N_x}=L\quad\text{and}\quad
 h_x=\max_{1\le j\le N_x}(x_j-x_{j-1}),
\]
and define the conforming cubic Hermite space
\[
 V_{h_x}=\left\{w_{h_x}\in C^1([0,L]):
 \begin{aligned}
  &w_{h_x}|_{[x_{j-1},x_j]}\in\mathbb P_3
     \quad(1\le j\le N_x),\\
  &w_{h_x}(0)=w_{h_x}'(0)=0
 \end{aligned}\right\}.
\]
Then $V_{h_x}\subset H_E^2(0,L)$. Standard conforming interpolation estimates for these elements are reviewed in \cite{Ciarlet2002}. The same beam family was used in
\cite{AnasseBomissoKpata2026}.
For the history variable, let
\[
 0=s_0<s_1<\cdots<s_N=S\quad\text{and}\quad h=\max_j(s_j-s_{j-1}),
\]
and set
\[
 V_h^s=\{\eta_h\in C([0,S]):
 \eta_h|_{[s_{j-1},s_j]}\in\mathbb P_1,
 \ \eta_h(0)=0\}.
\]
We use the weighted inner product
\[
 (u,v)_{m,S}=\int_0^S m(s)u(s)\overline{v(s)}\,ds
\]
and denote by $P_h^m:\Ms\to V_h^s$ the corresponding orthogonal projection. The discrete transport operator is
\[
 \mathcal G_{S,h}\eta_h=P_h^m(-\eta_h'),\quad\text{for all }\eta_h\in V_h^s.
\]

The finite-dimensional augmented space is
\[
 \Xks^{h_x,h}=V_{h_x}\times V_{h_x}\times V_h^s
\]
with the energy inner product inherited from \eqref{eq:inner-product}, with the history integral restricted to $(0,S)$. Define $\Aks^{h_x,h}(w_{h_x},v_{h_x},z_h)=(v_{h_x},a_{h_x},g_h)$, where $a_{h_x}\in V_{h_x}$ and $g_h\in V_h^s$ are uniquely determined by
\begin{align*}
 (a_{h_x},\phi_{h_x})_{L^2}
 ={}&-(w_{h_x}'',\phi_{h_x}'')_{L^2}
 -\zeta(v_{h_x},\phi_{h_x})_{L^2}
 -\delta^{(3)}w_{h_x}'(L)\overline{\phi_{h_x}'(L)}\\
 &-\delta^{(1)}v_{h_x}(L)\overline{\phi_{h_x}'(L)}
 -\delta^{(2)}v_{h_x}'(L)\overline{\phi_{h_x}'(L)}\\
 &-\gamma\bigl(v_{h_x}(L)+J_Sz_h\bigr)\overline{\phi_{h_x}(L)},\quad\text{for all }\phi_{h_x}\in V_{h_x},
\end{align*}
and
\[
 (g_h,\eta_h)_{m,S}
 =\bigl(v_{h_x}(L)-z_h',\eta_h\bigr)_{m,S},\quad\text{for all }\eta_h\in V_h^s.
\]
The two mass forms are positive definite, so the variational equations
uniquely define a finite-dimensional operator. The weighted integrals defining the scheme are taken exact. With numerical quadrature, the identities below hold up to the quadrature error. The following identity
holds on the kernel of the discrete memory moment:
\begin{equation}
 J_Sz_h=0
 \quad\Longrightarrow\quad
 \Aks^{h_x,h}(0,0,z_h)
 =(0,0,\mathcal G_{S,h}z_h).
 \label{eq:fe-pure-history}
\end{equation}

\begin{lemma}[Dissipativity of the Galerkin generator]
\label{lem:fe-dissipative}
Under \eqref{eq:structural}, the finite-dimensional operator
$\Aks^{h_x,h}$ is dissipative. For $U_h=(w_{h_x},v_{h_x},z_h)$,
\begin{equation}\label{eq:fe-dissipation}
 \begin{aligned}
 \Real\langle\Aks^{h_x,h}U_h,U_h\rangle
 ={}&-\zeta\|v_{h_x}\|_{L^2}^2
 -\mathbf v_{L,h}^*Q_L\mathbf v_{L,h}\\
 &-\frac{\gamma}{2}m(S)|z_h(S)|^2
 +\frac{\gamma}{2}\int_0^S m'(s)|z_h(s)|^2\,ds\le0,
 \end{aligned}
\end{equation}
where $\mathbf v_{L,h}=(v_{h_x}(L),v_{h_x}'(L))^\top$.
Consequently $\|e^{t\Aks^{h_x,h}}\|\le1$ for every $t\ge0$.
\end{lemma}

\begin{proof}
Take $\phi_{h_x}=v_{h_x}$ in the variational definition of $a_{h_x}$ and
$\eta_h=z_h$ in that of $g_h$. The stiffness terms cancel in real part
against the beam component of the energy inner product, as do the
$\delta^{(3)}$ terms. The mixed memory terms are
$-\gamma\Real(J_Sz_h\overline{v_{h_x}(L)})$ and
$\gamma\Real(v_{h_x}(L)\overline{J_Sz_h})$, which cancel.
The remaining physical boundary quadratic form is
$-\mathbf v_{L,h}^*Q_L\mathbf v_{L,h}$.
Since $z_h\in H^1(0,S)$ and $z_h(0)=0$,
\[
 -\Real(z_h',z_h)_{m,S}
 =-\frac12m(S)|z_h(S)|^2
 +\frac12\int_0^Sm'(s)|z_h(s)|^2\,ds.
\]
This gives \eqref{eq:fe-dissipation}. Positivity of $Q_L$ and
$m'\le0$ prove dissipativity. In finite dimension the operator
generates a semigroup, which is contractive by the energy identity.
\end{proof}

\begin{lemma}[Uniform weighted interpolation of the corrected packets]\label{lem:fe-interpolation}
Fix $\beta\in\R$. Let $z_S=z_{S,0}$ be given by Lemma~\ref{lem:moment-cancellation} if $\beta=0$, and $z_S=z_{S,\beta}$ be given by \eqref{eq:corrected-packet} if $\beta\ne0$. There exist $S_0>0$, $h_0>0$, and $C_\beta>0$, independent of $S\ge S_0$ and $0<h\le h_0$, such that the nodal piecewise-linear interpolant $I_hz_S$ satisfies
\begin{equation*}
 \|I_hz_S-z_S\|_{\Ms}
 +h\|(I_hz_S-z_S)'\|_{\Ms}
 \le C_\beta h^2.
\end{equation*}
The same estimate holds for the fixed correction function $\psi$.
\end{lemma}

\begin{proof}
For $\beta\ne0$, the corrected packet is a remote oscillatory packet plus a multiple $O_\beta(S^{-1-p/2})$ of $\psi$. Differentiating twice and using the normalization, the weight comparability on $I_S$ and $\ell_S=S/8$ gives
\begin{equation}
 \sup_{S\ge S_0}
 \bigl(\|z_S\|_{\Ms}+\|z_S'\|_{\Ms}+\|z_S''\|_{\Ms}\bigr)
 \le C_\beta.
 \label{eq:uniform-packet-h2}
\end{equation}
At $\beta=0$ the same bound follows from \eqref{eq:zero-moment-direct}. For a history cell $K=[s_{j-1},s_j]$ of length $h_K=s_j-s_{j-1}$, let $I_K$ be the $\mathbb P_1$ Lagrange interpolant on $K$ and let $L^2_m(K)$ be the $L^2$ space on $K$ with weight $m$. On each cell with $h_K\le h_0$, $|m'/m|\le p+1$ gives $\sup_Km/\inf_Km\le e^{(p+1)h_0}$. The unweighted estimates $\|z-I_Kz\|_{L^2(K)}\le Ch_K^2\|z''\|_{L^2(K)}$ and $\|(z-I_Kz)'\|_{L^2(K)}\le Ch_K\|z''\|_{L^2(K)}$ \cite{Ciarlet2002} therefore hold in $L^2_m(K)$ with an extra factor $e^{(p+1)h_0/2}$. Summing over cells gives $\|I_hz-z\|_{\Ms}+h\|(I_hz-z)'\|_{\Ms}\le Ch^2\|z''\|_{\Ms}$ with $C$ independent of $S$, and \eqref{eq:uniform-packet-h2} concludes.
\end{proof}

\begin{theorem}[Finite-element persistence of the pseudospectral bounds]\label{thm:fe-pseudospectral-transfer}
Fix $\beta\in\R$. There exist $S_0>0$, $h_0>0$, and $C_\beta>0$, independent of the truncation length $S$, the cubic-Hermite beam mesh size $h_x$, and the history mesh size $h$, such that for every $S\ge S_0$, every cubic-Hermite beam mesh of size $h_x$, and every $0<h\le h_0$, there is a vector $U_{S,h_x,h}\in\Xks^{h_x,h}$ satisfying
\begin{equation}
 \|U_{S,h_x,h}\|_{\Xks^{h_x,h}}=1
 \quad\text{and}\quad
 \|(i\beta I-\Aks^{h_x,h})U_{S,h_x,h}\|_{\Xks^{h_x,h}}
 \le C_\beta\bigl(S^{-1}+h\bigr).
 \label{eq:fe-pseudospectral-bound}
\end{equation}
Consequently,
\[
 i\beta\in
 \sigma_{C_\beta(S^{-1}+h)}(\Aks^{h_x,h}).
\]
The states have zero beam components, so the construction and its constants do not depend on the beam mesh $h_x$. At zero frequency, Proposition~\ref{prop:fe-zero-frequency} removes the term $h$.
\end{theorem}

\begin{proof}
Let $I_h$ denote the nodal interpolant into $V_h^s$, and put $\psi_h=I_h\psi$. Since $S>2$ and $\psi$ is supported in $(1,2)$,
\[
 J_S\psi=J\psi>0.
\]
By Lemma~\ref{lem:fe-interpolation} and the continuity of $J_S$,
\[
 |J_S(\psi_h-\psi)|\le
 \left(\int_0^S m(s)\,ds\right)^{1/2}
 \|\psi_h-\psi\|_{\Ms}
 \le Ch^2.
\]
Hence, after decreasing $h_0$ if necessary,
\[
 J_S\psi_h\ge\frac12J\psi>0
\]
uniformly in $S$. Since $J_Sz_S=0$, we also have
\[
 |J_SI_hz_S|
 =|J_S(I_hz_S-z_S)|\le Ch^2.
\]
Nodal interpolation does not preserve the moment, even at $\beta=0$. We define the exactly constrained interpolant
\[
 \widehat z_{S,h}
 =I_hz_S-
 \frac{J_SI_hz_S}{J_S\psi_h}\,\psi_h.
\]
Then $J_S\widehat z_{S,h}=0$, and Lemma~\ref{lem:fe-interpolation} gives
\begin{equation}
 \|\widehat z_{S,h}-z_S\|_{\Ms}
 +\|(\widehat z_{S,h}-z_S)'\|_{\Ms}
 \le C_\beta h.
 \label{eq:fe-graph-approx}
\end{equation}
In particular, $\|\widehat z_{S,h}\|_{\Ms}\to1$ uniformly in $S$, so for small $h$ we may set
\[
 z_{S,h}=
 \frac{\widehat z_{S,h}}{\|\widehat z_{S,h}\|_{\Ms}}.
\]
The normalization changes \eqref{eq:fe-graph-approx} only by another $O(h)$ term.

Because $P_h^m$ is an orthogonal projection in $\Ms$ and $z_{S,h}\in V_h^s$,
\[
 (i\beta I-\mathcal G_{S,h})z_{S,h}
 =P_h^m\bigl(i\beta z_{S,h}+z_{S,h}'\bigr).
\]
Therefore, using the contractivity of $P_h^m$, Proposition~\ref{prop:finite-history-pseudo}, and \eqref{eq:fe-graph-approx},
\[
 \begin{aligned}
 \|(i\beta I-\mathcal G_{S,h})z_{S,h}\|_{\Ms}
 &\le \|(i\beta I-\mathcal G_S)z_S\|_{\Ms}\\
 &\quad+|\beta|\|z_{S,h}-z_S\|_{\Ms}
 +\|(z_{S,h}-z_S)'\|_{\Ms}\\
 &\le C_\beta\bigl(S^{-1}+h\bigr).
 \end{aligned}
\]
Finally define
\[
 U_{S,h_x,h}=\gamma^{-1/2}(0,0,z_{S,h}).
\]
Its discrete energy norm is one. Since $J_Sz_{S,h}=0$, the exact identity \eqref{eq:fe-pure-history} yields
\[
 \Aks^{h_x,h}U_{S,h_x,h}
 =\gamma^{-1/2}(0,0,\mathcal G_{S,h}z_{S,h}).
\]
Thus \eqref{eq:fe-pseudospectral-bound} follows. The pseudospectral inclusion is immediate from the definition.
\end{proof}

\begin{proposition}[Sharper discrete residual at zero frequency]
\label{prop:fe-zero-frequency}
Let $p>1$ and assume the polynomial memory kernel and its weight defined in \eqref{eq:kernel}.
There are $S_0>0$, $h_0\in(0,1]$ and $C>0$, independent of
$S\ge S_0$, the beam mesh $h_x$ and the history mesh $0<h\le h_0$,
such that a unit vector $U_{S,h_x,h}^{(0)}\in\Xks^{h_x,h}$ satisfies
\begin{equation}\label{eq:fe-zero-frequency-bound}
 \|\Aks^{h_x,h}U_{S,h_x,h}^{(0)}\|_{\Xks^{h_x,h}}
 \le C\bigl(S^{-1}+hS^{-2}\bigr)\le CS^{-1}.
\end{equation}
The vector has vanishing beam components.
\end{proposition}

\begin{proof}
Take $z_S=z_{S,0}$ from Lemma~\ref{lem:moment-cancellation}, so that $J_Sz_S=0$, $\|z_S'\|_{\Ms}\le CS^{-1}$ and $\|z_S''\|_{\Ms}\le CS^{-2}$. The cellwise estimate in the proof of Lemma~\ref{lem:fe-interpolation} gives $\|I_hz_S-z_S\|_{\Ms}\le Ch^2S^{-2}$ and $\|(I_hz_S-z_S)'\|_{\Ms}\le ChS^{-2}$. With $\psi_h=I_h\psi$ as in Theorem~\ref{thm:fe-pseudospectral-transfer}, the moment-corrected interpolant $\widehat z_{S,h}=I_hz_S-\alpha_{S,h}\psi_h$, $\alpha_{S,h}=J_SI_hz_S/J_S\psi_h=O(h^2S^{-2})$, satisfies $J_S\widehat z_{S,h}=0$ and the same two bounds. After normalization, $\|z_{S,h}'\|_{\Ms}\le C(S^{-1}+hS^{-2})$. Since $P_h^m$ is contractive, $\|\mathcal G_{S,h}z_{S,h}\|_{\Ms}\le\|z_{S,h}'\|_{\Ms}$, and $U_{S,h_x,h}^{(0)}=\gamma^{-1/2}(0,0,z_{S,h})$ satisfies \eqref{eq:fe-zero-frequency-bound} by \eqref{eq:fe-pure-history}.
\end{proof}

\begin{corollary}[Failure of uniform exponential stability after approximation]
\label{cor:fe-nonuniform}
Assume \eqref{eq:structural}. Let $S_0$ and $h_0$ be the thresholds of Theorem~\ref{thm:fe-pseudospectral-transfer} and Proposition~\ref{prop:fe-zero-frequency}, and let $S\ge S_0$ and $0<h\le h_0$. Suppose the semigroup of a fixed
$\Aks^{h_x,h}$ satisfies $\|e^{t\Aks^{h_x,h}}\|\le Me^{-\omega t}$
for $t\ge0$, with $M\ge1$ and $\omega>0$. Then for every fixed
$\beta\in\R$,
\begin{equation}\label{eq:rate-bound}
 \frac{\omega}{M}\le C_\beta(S^{-1}+h).
\end{equation}
Moreover, at zero frequency the stronger bound
\begin{equation}\label{eq:rate-bound-zero}
 \frac{\omega}{M}\le CS^{-1}
\end{equation}
holds with $C$ independent of $S\ge S_0$, of the beam mesh $h_x$ and of $0<h\le h_0$. Consequently no common $M\ge1$ and $\omega>0$
can give an exponential bound as $S\to\infty$, even if the
history step $h\in(0,h_0]$ is kept fixed.
\end{corollary}

\begin{proof}
The assumed exponential estimate gives
$i\R\subset\rho(\Aks^{h_x,h})$ and, by the Laplace integral,
\[
 \|(i\beta I-\Aks^{h_x,h})^{-1}\|
 \le\int_0^\infty Me^{-\omega t}\,dt=M/\omega.
\]
Take the unit vector $U_{S,h_x,h}$ from
Theorem~\ref{thm:fe-pseudospectral-transfer}. Then
\[
 1\le\frac M\omega
 \|(i\beta I-\Aks^{h_x,h})U_{S,h_x,h}\|
 \le\frac M\omega C_\beta(S^{-1}+h),
\]
which proves \eqref{eq:rate-bound}. At $\beta=0$, replace that
vector by $U_{S,h_x,h}^{(0)}$ from
Proposition~\ref{prop:fe-zero-frequency}. The same resolvent
argument and \eqref{eq:fe-zero-frequency-bound} yield
\eqref{eq:rate-bound-zero}. If $M$ and $\omega$ were
independent of $S$, this would be impossible as $S\to\infty$,
including along families with a fixed $0<h\le h_0$.
\end{proof}

Because $\Aks^{h_x,h}$ is not normal, these residual bounds constrain the ratio $\omega/M$ but do not locate its eigenvalues.

\subsection{Zero-frequency resolvent growth and optimal time scale}
\label{subsec:zero-frequency-growth}

The zero-frequency residuals of Proposition~\ref{prop:fe-zero-frequency} lead to a quantitative statement for the full coupled generator: its resolvent at zero grows exactly linearly in $S$, both before and after discretization. The matching upper bound is uniform in the two meshes. Write
\[
 a_0(w,\phi)=\int_0^L w''\overline{\phi''}\,dx+\delta^{(3)}w'(L)\overline{\phi'(L)},
\]
so that the beam part of the energy norm is $\|w\|_{a_0}^2=a_0(w,w)$, and let $\mathcal G_{S,h}=P_h^m(-\partial_s)$ on $V_h^s$.

\begin{theorem}[Linear zero-frequency resolvent growth]
\label{thm:zero-frequency-growth}
Assume \eqref{eq:structural}. There are $S_0>2$, $h_0\in(0,1]$ and constants $0<c\le C$, independent of $S$, of the beam mesh $h_x$ and of $0<h\le h_0$, such that $0\in\rho(\Aks)\cap\rho(\Aks^{h_x,h})$ and, for $S\ge S_0$,
\begin{equation}
 cS\le\|\Aks^{-1}\|_{\mathcal L(\Xks)}\le CS
 \quad\text{and}\quad
 cS\le\|(\Aks^{h_x,h})^{-1}\|_{\mathcal L(\Xks^{h_x,h})}\le CS.
 \label{eq:zero-frequency-growth}
\end{equation}
\end{theorem}

\begin{proof}
\emph{Transport coercivity.} For $z\in D(\mathcal G_S)$, the condition $z(0)=0$ and the identity $-m'(s)=\frac{p+1}{1+s}m(s)$ give
\[
 \Real(z',z)_{m,S}
 =\tfrac12m(S)|z(S)|^2+\tfrac12\int_0^S[-m'(s)]|z(s)|^2\,ds
 \ge\frac{p+1}{2(1+S)}\|z\|_{\Ms}^2.
\]
The same computation applies to every $z_h\in V_h^s$, and $\Real(\mathcal G_{S,h}z_h,z_h)_{m,S}=-\Real(z_h',z_h)_{m,S}$ because $P_h^m$ is the orthogonal projection. Hence $-\mathcal G_S$ and $-\mathcal G_{S,h}$ are invertible, with
\begin{equation}
 \|(-\mathcal G_S)^{-1}\|,\ \|(-\mathcal G_{S,h})^{-1}\|\le\frac{2(1+S)}{p+1}.
 \label{eq:transport-inverse}
\end{equation}
For $\mathcal G_S$, the formula $z(s)=\int_0^sf(r)\,dr$ proves surjectivity. In the finite-dimensional space $V_h^s$, injectivity suffices.

\emph{Upper bound.} Let $F=(f,g,\vartheta)$ and solve $-\Aks U=F$ with $U=(w,v,z)$. The first component gives $v=-f$. The third gives $(-\mathcal G_S)z=\vartheta+v(L)$, where $v(L)$ is the constant history. Since $\int_0^Sm\le b$ and $|f(L)|\le C\|f\|_{a_0}$, \eqref{eq:transport-inverse} yields $\|z\|_{\Ms}\le C(1+S)\|F\|$. The beam component solves the static problem
\[
 a_0(w,\phi)=(g,\phi)_{L^2}-\zeta(v,\phi)_{L^2}
 -\delta^{(1)}v(L)\overline{\phi'(L)}
 -\delta^{(2)}v'(L)\overline{\phi'(L)}
 -\gamma\bigl(v(L)+J_Sz\bigr)\overline{\phi(L)}
\]
for all $\phi\in H_E^2(0,L)$. By \eqref{eq:Jbound} restricted to $(0,S)$, $|J_Sz|\le\sqrt b\|z\|_{\Ms}$. The Lax--Milgram lemma gives $\|w\|_{a_0}\le C(1+S)\|F\|$, and one-dimensional elliptic regularity gives $U\in D(\Aks)$. The discrete problem has the same structure: $v_{h_x}=-f_{h_x}$, $(-\mathcal G_{S,h})z_h=\vartheta_h+P_h^m(v_{h_x}(L))$, and the static beam problem is posed on $V_{h_x}$, where $a_0$ has the same coercivity constant. This proves the upper bounds in \eqref{eq:zero-frequency-growth}.

\emph{Lower bound.} The function $z_{S,0}$ of Lemma~\ref{lem:moment-cancellation} satisfies $J_Sz_{S,0}=0$ and $\|z_{S,0}'\|_{\Ms}\le CS^{-1}$. The vector $U=\gamma^{-1/2}(0,0,z_{S,0})$ belongs to $D(\Aks)$, has norm one, and $\|\Aks U\|\le CS^{-1}$. Hence $\|\Aks^{-1}\|\ge S/C$. Proposition~\ref{prop:fe-zero-frequency} gives the discrete lower bound in the same way.
\end{proof}

\begin{remark}[The restriction $h\le h_0$]
\label{rem:coarse-history-mesh}
The lower bound in \eqref{eq:zero-frequency-growth} fails for arbitrarily coarse history meshes. Let $p=3$, $b=1$, and take the single cell $[0,S]$, so that $V_h^s=\operatorname{span}\{e_S\}$ with $e_S(s)=s/S$. Put
\[
 A_j(S)=\int_0^Ss^jm(s)\,ds,\quad\text{for }j=1,2.
\]
The weighted projection gives
\[
 -\mathcal G_{S,h}e_S=\frac{A_1(S)}{A_2(S)}\,e_S
 \quad\text{and}\quad
 \|P_h^m1\|_{\Ms}=\frac{A_1(S)}{\sqrt{A_2(S)}}.
\]
Since $A_1(S)\to1/2$ and $A_2(S)\to1$, both $\|(-\mathcal G_{S,h})^{-1}\|$ and $\|P_h^m1\|_{\Ms}$ stay bounded. The argument of the upper bound then gives $\sup_{S\ge S_0}\|(\Aks^{h_x,h})^{-1}\|<\infty$. Thus the condition $h\le h_0$ in Theorem~\ref{thm:zero-frequency-growth} is necessary.
\end{remark}

\begin{corollary}[Plateau of the semigroup norms]
\label{cor:plateau}
Under the assumptions of Theorem~\ref{thm:zero-frequency-growth}, there is $C_*>0$, independent of $S$, $h_x$ and $0<h\le h_0$, such that
\[
 \|e^{t\Aks}\|\ge1-\frac{C_*t}{S}\quad\text{and}\quad
 \|e^{t\Aks^{h_x,h}}\|\ge1-\frac{C_*t}{S},\quad\text{for all }t\ge0.
\]
In particular, for $0<\varepsilon<1$, neither semigroup norm falls below $1-\varepsilon$ before the time $\varepsilon S/C_*$.
\end{corollary}

\begin{proof}
Let $U$ be the unit vector used in the lower bound above. Since the semigroups are contractive,
$\|e^{t\Aks}U-U\|\le\int_0^t\|e^{r\Aks}\Aks U\|\,dr\le t\|\Aks U\|\le C_*t/S$, and the same argument applies to $\Aks^{h_x,h}$.
\end{proof}

For every admissible mesh and large $S$, some unit-energy states retain a fixed fraction of their energy for times of order $S$. Their mechanical components vanish, so this plateau is produced by the prehistory alone.

With strictly positive distributed and boundary damping, a Galerkin-compatible Lyapunov functional gives the matching upper bound \cite{KoffiAnasseDansouBomisso2026}.

\begin{theorem}[Optimal uniform exponential time scale]
\label{thm:optimal-time-scale}
Assume \eqref{eq:structural}, $\zeta>0$ and that $Q_L$ is positive definite. There are $M\ge1$ and $c>0$, independent of $S>2$, $h_x$ and $h$, such that
\begin{equation}
 \|e^{t\Aks}\|\le Me^{-ct/(1+S)}\quad\text{and}\quad
 \|e^{t\Aks^{h_x,h}}\|\le Me^{-ct/(1+S)},\quad\text{for all }t\ge0.
 \label{eq:uniform-time-scale}
\end{equation}
Moreover, for every $M_*\ge M$, the largest rate $\omega_S(M_*)$ for which $\|e^{t\Aks}\|\le M_*e^{-\omega t}$ holds for all $t\ge0$ satisfies, for $S\ge S_0$,
\begin{equation}
 \frac{c}{1+S}\le\omega_S(M_*)\le\frac{2C_*\log(2M_*)}{S},
 \label{eq:optimal-rate}
\end{equation}
where $C_*$ is the constant of Corollary~\ref{cor:plateau}. The same bounds hold for $\Aks^{h_x,h}$ uniformly in $h_x$ and $0<h\le h_0$.
\end{theorem}

\begin{proof}
We first consider classical solutions. The Galerkin proof uses the admissible tests $\phi_{h_x}=w_{h_x}$ and $\eta_h=z_h$. The continuous result extends by density. Set $B^2=a_0(w,w)$, $V=\|v\|_{L^2}$, $Z=\|z\|_{\Ms}$ and $\mathbf v_L=(v(L),v'(L))^\top$. Then $E=\frac12(B^2+V^2+\gamma Z^2)$. Let $\lambda_L>0$ be the smallest eigenvalue of $Q_L$.

\emph{Dissipation.} By the energy identity and $-m'(s)=\frac{p+1}{1+s}m(s)\ge\frac{p+1}{1+S}m(s)$ on $(0,S)$,
\begin{equation}
 E'\le-\zeta V^2-\lambda_L|\mathbf v_L|^2-\frac{\gamma d_0}{1+S}Z^2,\quad\text{with}\quad d_0=\frac{p+1}{2}.
 \label{eq:time-scale-dissipation}
\end{equation}

\emph{Beam multiplier.} Let $K=\Real(w,v)_{L^2}$. Testing the beam equation with $w$ gives
\[
 K'=V^2-B^2-\zeta\Real(v,w)_{L^2}
 -\delta^{(1)}\Real\bigl(v(L)\overline{w'(L)}\bigr)
 -\delta^{(2)}\Real\bigl(v'(L)\overline{w'(L)}\bigr)
 -\gamma\Real\bigl((v(L)+J_Sz)\overline{w(L)}\bigr).
\]
On the fixed interval $(0,L)$, $\|w\|_{L^2}\le C_PB$ and $|w(L)|+|w'(L)|\le C_TB$. Since $|J_Sz|\le\sqrt b\,Z$, Young's inequality gives $C_0>0$, independent of $S$ and of the meshes, such that
\[
 K'\le-\tfrac12B^2+C_0\bigl(V^2+|\mathbf v_L|^2+\gamma Z^2\bigr)\quad\text{and}\quad |K|\le C_PE.
\]

\emph{Lyapunov functional.} Choose $\eta>0$ with $\eta C_P\le\frac12$ and $\eta C_0\le\frac12\min\{\zeta,\lambda_L,d_0\}$, and set $F=E+\frac{\eta}{1+S}K$. Then $\frac12E\le F\le\frac32E$ and
\[
 F'\le-\frac{\eta}{2(1+S)}B^2-\frac\zeta2V^2-\frac{\lambda_L}{2}|\mathbf v_L|^2-\frac{\gamma d_0}{2(1+S)}Z^2
 \le-\frac{d}{1+S}E\le-\frac{2d}{3(1+S)}F,
\]
with $d=\min\{\eta,\zeta,d_0\}$. Hence $E(t)\le3e^{-2dt/(3(1+S))}E(0)$, which gives \eqref{eq:uniform-time-scale} with $M=\sqrt3$ and $c=d/3$, and the lower bound in \eqref{eq:optimal-rate}.

\emph{Optimality.} By Corollary~\ref{cor:plateau}, both semigroup norms are at least $\frac12$ at $t_S=S/(2C_*)$. A bound $M_*e^{-\omega t}$ therefore requires $\frac12\le M_*e^{-\omega t_S}$, that is, $\omega\le 2C_*\log(2M_*)/S$.
\end{proof}

The rate describes the optimal uniform time scale, not the spectral abscissa of each fixed generator. The optimality concerns the norm on the augmented space. The vectors that saturate it have zero displacement and zero velocity.

For the continuous generator, a boundary multiplier removes the need for distributed damping \cite{LazzariNibbi2012}. The constants remain uniform in $S$.

\begin{theorem}[Optimal time scale without distributed damping]
\label{thm:time-scale-zeta-zero}
Assume \eqref{eq:structural}, that $Q_L$ is positive definite, and that $\delta^{(3)}L\le2$. There are $M\ge1$ and $c>0$, independent of $S>2$, such that
\[
 \|e^{t\Aks}\|\le Me^{-ct/(1+S)},\quad\text{for all }t\ge0.
\]
Consequently the bounds \eqref{eq:optimal-rate} hold for $\Aks$ under these assumptions, in particular for $\zeta=0$.
\end{theorem}

\begin{proof}
Retain the notation of Theorem~\ref{thm:optimal-time-scale}. Set $a=w'(L)$, $c=w''(L)$ and $r=\delta^{(1)}v(L)+\delta^{(2)}v'(L)$. The boundary conditions give $c=-\delta^{(3)}a-r$ and $w'''(L)=\gamma(v(L)+J_Sz)$. Since $w'(0)=0$, we have $|a|^2\le L\|w''\|_{L^2}^2$ and $\|w'\|_{L^2}^2\le\frac{L^2}{2}\|w''\|_{L^2}^2$. Moreover $B^2\le(1+\delta^{(3)}L)\|w''\|_{L^2}^2\le3\|w''\|_{L^2}^2$.

\emph{Boundary multiplier.} For classical solutions let $P=\Real(xw',v)_{L^2}$. Then $|P|\le C_PE$. Using $v_t=-w^{(4)}-\zeta v$, $v(0)=0$, $w(0)=w'(0)=0$ and integration by parts,
\[
 \Real(xv',v)_{L^2}=\tfrac L2|v(L)|^2-\tfrac12V^2,
\]
\[
 \Real(xw',-w^{(4)})_{L^2}
 =-L\Real\bigl(a\overline{w'''(L)}\bigr)+\Real(a\overline c)
 -\tfrac32\|w''\|_{L^2}^2+\tfrac L2|c|^2.
\]
Inserting the boundary conditions,
\[
 \tfrac L2|c|^2+\Real(a\overline c)
 =\delta^{(3)}\Bigl(\tfrac{\delta^{(3)}L}{2}-1\Bigr)|a|^2
 +(\delta^{(3)}L-1)\Real(a\overline r)+\tfrac L2|r|^2,
\]
and the first term is nonpositive because $\delta^{(3)}L\le2$.
To control the remaining mixed boundary terms, set
\[
 d_\delta=\bigl((\delta^{(1)})^2+(\delta^{(2)})^2\bigr)^{1/2},\quad
 A_0=|\delta^{(3)}L-1|d_\delta+\gamma L\quad\text{and}\quad
 B_0=\gamma L\sqrt b.
\]
Since $|r|\le d_\delta|\mathbf v_L|$, $|a|\le\sqrt L\,\|w''\|_{L^2}$ and $|J_Sz|\le\sqrt b\,Z$, Young's inequality yields
\begin{align*}
 &\bigl|(\delta^{(3)}L-1)\Real(a\overline r)
     -\gamma L\Real\bigl(a\overline{v(L)+J_Sz}\bigr)\bigr|\\
 &\qquad\le\sqrt L\,\|w''\|_{L^2}
       (A_0|\mathbf v_L|+B_0Z)\\
 &\qquad\le\tfrac12\|w''\|_{L^2}^2
       +L A_0^2|\mathbf v_L|^2+L B_0^2 Z^2.
\end{align*}
The terms $\frac L2|v(L)|^2$ and $\frac L2|r|^2$ are bounded by
$\frac L2(1+d_\delta^2)|\mathbf v_L|^2$. Moreover,
$\|xw'\|_{L^2}\le L^2\|w''\|_{L^2}/\sqrt2$, so
\[
 \zeta|\Real(xw',v)_{L^2}|
 \le\tfrac14\|w''\|_{L^2}^2+\tfrac{\zeta^2L^4}{2}V^2.
\]
Consequently the coefficient of $\|w''\|_{L^2}^2$ in $P'$ is at most
$-3/2+1/2+1/4=-3/4$. As $B^2\le3\|w''\|_{L^2}^2$, there are constants
$C_1,C_2>0$, independent of $S$, such that
\[
 P'\le-\tfrac12V^2-\tfrac1{12}B^2
 +C_1\bigl(|\mathbf v_L|^2+\gamma Z^2\bigr)+C_2\zeta^2V^2.
\]

\emph{Lyapunov functional.} Choose $\eta>0$ with $\eta C_P\le\frac12$, $\eta C_1\le\frac12\min\{\lambda_L,d_0\}$ and $\eta C_2\zeta\le\frac12$, and set $F=E+\frac{\eta}{1+S}P$. Then $\frac12E\le F\le\frac32E$, and \eqref{eq:time-scale-dissipation}, which holds for every $\zeta\ge0$, gives
\[
 F'\le-\frac{\eta}{1+S}\Bigl(\tfrac12V^2+\tfrac1{12}B^2\Bigr)-\frac{\gamma d_0}{2(1+S)}Z^2
 \le-\frac{d}{1+S}E,\quad\text{with}\quad d=\min\{\eta/6,\,d_0\}.
\]
Gronwall's lemma concludes as before, and the upper bound in \eqref{eq:optimal-rate} follows from Corollary~\ref{cor:plateau}.
\end{proof}

The multiplier $xw'$ does not belong to the Hermite space, so Theorem~\ref{thm:time-scale-zeta-zero} concerns the continuous truncated generator. The degenerate boundary-damping case $\delta^{(1)}=\delta^{(2)}=0$ is stated explicitly below.

\begin{corollary}[Pure tip-velocity feedback]
\label{cor:tip-velocity-feedback}
Assume $p>1$, $b,\gamma>0$, $\delta^{(1)}=\delta^{(2)}=0$, $\delta^{(3)}\ge0$ and $\delta^{(3)}L\le2$. For every fixed $\zeta\ge0$, the continuous truncated generator satisfies
\[
 \|e^{t\Aks}\|\le M e^{-ct/(1+S)},\quad\text{for all }t\ge0,
\]
with $M,c>0$ independent of $S>2$. If moreover $\zeta>0$, the same estimate holds for $\Aks^{h_x,h}$ uniformly in the two meshes. In each applicable case the order $(1+S)^{-1}$ is optimal for estimates with a prefactor independent of $S$.
\end{corollary}
\begin{proof}
For $\delta^{(1)}=\delta^{(2)}=0$, the dissipation identity controls $\gamma|v(L)|^2$ rather than the full trace vector. In the proof of Theorem~\ref{thm:time-scale-zeta-zero} we have $r=0$. The explicit bound on the boundary multiplier above then contains only $|v(L)|^2$ and $Z^2$, both absorbed by the dissipation after choosing the Lyapunov coefficient sufficiently small. This proves the continuous estimate for $\zeta\ge0$. When $\zeta>0$, the multiplier $K=\Real(w,v)_{L^2}$ used in Theorem~\ref{thm:optimal-time-scale} likewise produces only the trace $v(L)$, so its argument applies without change in the continuous and Galerkin spaces. The lower bound follows from Corollary~\ref{cor:plateau}.
\end{proof}

\subsection{Numerical verification}
\label{subsec:numerical-illustration}

We use
\[
 p=2,\qquad b=1\quad\text{and}\quad \beta=1,
\]
with
\[
 \varphi(r)=
 \begin{cases}
  \exp\!\left(-\dfrac{1}{1-r^2}\right),& |r|<1,\\[1mm]
  0,& |r|\ge1,
 \end{cases}
 \quad\text{and}\quad
 \psi(s)=\varphi(2s-3),
\]
and correct the discrete moment using the same algebraic formula as in the proofs. In the computations, the moment coefficients are evaluated by ten-point Gauss--Legendre quadrature. The analytical construction uses the exact weighted moment. For the continuous remote packets, we factor out $m(R)$ and integrate the relative weight $m(R+\ell r)/m(R)$. This avoids applying a fixed absolute quadrature tolerance to integrals that become very small as $R$ increases.

At frequency $\beta=0$ we also checked the direct two-profile construction
of Lemma~\ref{lem:moment-cancellation} with $p=1.5$ and $b=1$.
Using the same smooth function $\varphi$ as above, the two profiles in
the rescaled variable $r$ are
\[
 \varphi_1(r)=\varphi\!\left(\frac{r+0.4}{0.25}\right)\quad\text{and}\quad
 \varphi_2(r)=\varphi\!\left(\frac{r-0.4}{0.25}\right).
\]
Their disjoint supports are $[-0.65,-0.15]$ and $[0.15,0.65]$. The coefficient $c_S$ cancels the weighted moment as in Lemma~\ref{lem:moment-cancellation}. The normalized residuals are:
\[
\begin{array}{c|ccccc}
S&64&128&256&512&1024\\\hline
\|\mathcal G_S z_{S,0}\|_{\Ms}
&0.87825&0.43914&0.21957&0.10979&0.05489\\
S\|\mathcal G_S z_{S,0}\|_{\Ms}
&56.2079&56.2093&56.2101&56.2104&56.2106
\end{array}
\]
The products approach a constant, confirming the predicted order $O(S^{-1})$ for $1<p<2$.

\begin{figure}[htbp]
 \centering
 \includegraphics[width=0.94\textwidth]{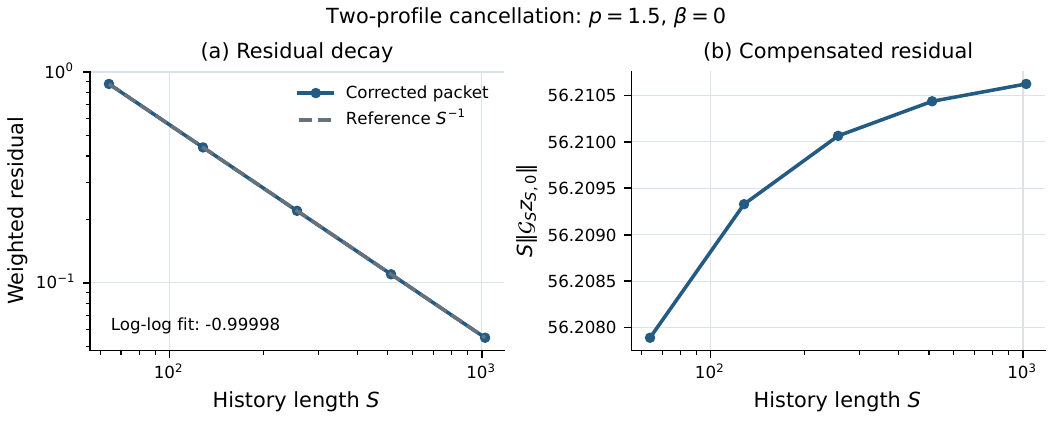}
 \caption{Two-profile moment cancellation at zero frequency for $p=1.5$ and $b=1$. Left: continuous translation residual and $S^{-1}$ reference. Right: compensated residual on a magnified scale.}
 \label{fig:zero-frequency-two-profile}
\end{figure}

For the finite-history packets, set
\[
 \varepsilon_S:=\|(i\beta I-\mathcal G_S)z_S\|_{\Ms}.
\]
With \(R_S=3S/4\), \(\ell_S=S/8\) and
\(S=32,64,128,256,512,1024\), the fitted slope of \(\varepsilon_S\) against \(S\) is \(-1.00181\). The products \(S\varepsilon_S\) stabilize near \(14.44\). This is consistent with the $O(S^{-1})$ estimate of Proposition~\ref{prop:finite-history-pseudo} at the fixed frequency $\beta=1$.

\begin{figure}[H]
 \centering
 \includegraphics[width=0.94\textwidth]{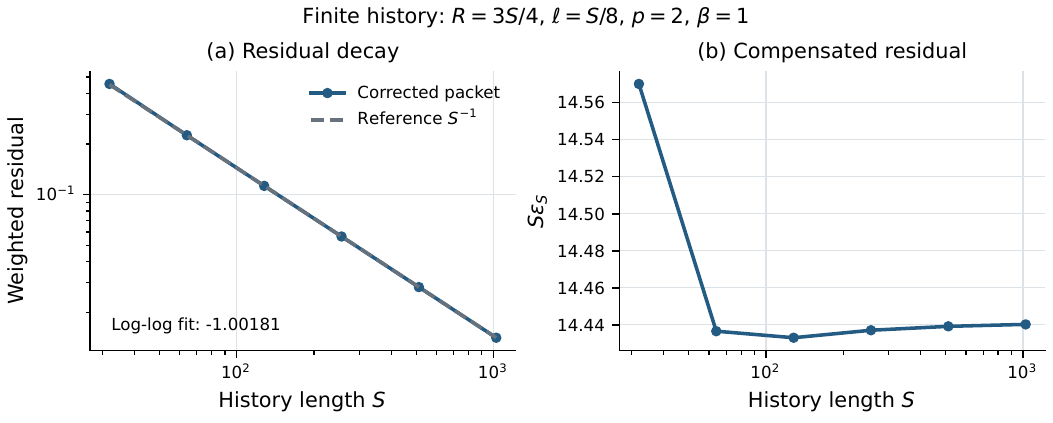}
 \caption{Finite-history packets for $p=2$, $b=1$, $\beta=1$. Left: residual and $S^{-1}$ reference. Right: compensated residual $S\varepsilon_S$ on a magnified vertical scale.}
 \label{fig:finite-history-pseudo}
\end{figure}

Finally, fix \(S=256\) and use uniform history meshes
\[
 h\in\left\{\frac12,\frac14,\frac18,\frac1{16},\frac1{32}\right\}.
\]
With ten-point Gauss--Legendre quadrature for the weighted \(\mathbb P_1\) matrices and the moment vector, the graph defect
\[
 D_h=\|z_{S,h}-z_S\|_{\Ms}
 +\|(z_{S,h}-z_S)'\|_{\Ms}
\]
has fitted slope \(1.02328\). Moreover,
\[
 \varepsilon_{S,h}
 =\|(i\beta I-\mathcal G_{S,h})z_{S,h}\|_{\Ms}
 \longrightarrow
 \varepsilon_S=0.05639489885.
\]
The graph defect is first-order, and the discrete residual converges to its finite-history counterpart.

\begin{figure}[H]
 \centering
 \includegraphics[width=0.48\textwidth]{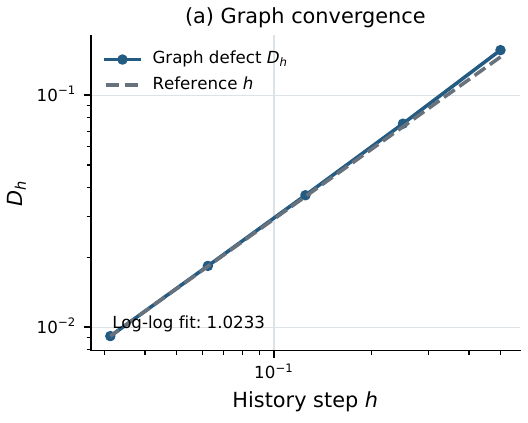}\hfill
 \includegraphics[width=0.48\textwidth]{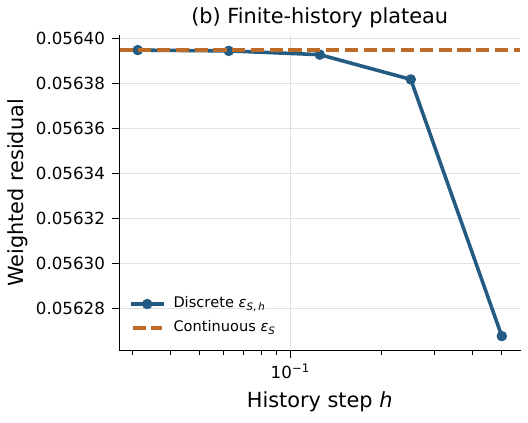}
 \caption{History-mesh verification for \(S=256\), \(p=2\), \(b=1\), and \(\beta=1\). Left: first-order graph defect. Right: convergence of the discrete residual to the finite-history plateau.}
 \label{fig:fe-history-mesh}
\end{figure}

The slopes agree with the residual estimates. We check the corrected discrete moment using higher-order quadrature on the same piecewise linear functions. For $S=256$ and $h\in\{1/2,1/4,1/8,1/16,1/32\}$,
the reference quadrature gives
$|J_S^{(48)}z_{S,h}|\le 1.5\times10^{-19}$, while
$|J_S^{(48)}z_{S,h}-J_S^{(10)}z_{S,h}|\le 1.5\times10^{-19}$.
The superscript gives the number of Gauss--Legendre points per element.
Identity~\eqref{eq:fe-pure-history} links these packets to the Galerkin generator.

\section{Conclusion}
\label{sec:conclusion}

Truncation turns remote-history packets into quantitative pseudospectral residuals, and the Hermite--$\mathbb P_1$ approximation retains the zero-frequency construction. The main quantitative results concern the truncated and discrete generators. Their zero-frequency resolvents have sharp growth of order $S$, uniformly in the meshes. Their semigroup norms retain a plateau of duration proportional to $S$.

Under strictly positive distributed and boundary damping, the matching exponential rate $1/(1+S)$ holds uniformly in both meshes. The continuous truncated beam has the same optimal time scale without distributed damping under Theorem~\ref{thm:time-scale-zeta-zero} or Corollary~\ref{cor:tip-velocity-feedback}. In that regime the mesh-uniform question remains open. The separate high-frequency Galerkin mechanism is conditional and discussed in Appendix~\ref{app:high-frequency}. The finite-history solutions converge strongly to their complete-history counterparts on bounded time intervals. The ancillary static transfer result and its fixed-$S$ Galerkin convergence are recorded in Appendix~\ref{app:static-transfer}.

\appendix
\section{Static memory-to-beam transfer at fixed history length}
\label{app:static-transfer}
The following static calculation is complementary to the mesh-uniform exponential time-scale theorem. It distinguishes the physical displacement response to a history source from the growth of the complete resolvent.

The vectors in Theorem~\ref{thm:zero-frequency-growth} have zero beam components. The beam response to a history forcing is described by the displacement block of the same inverse, which can be computed exactly. Let $\iota_S\vartheta=(0,0,\vartheta)$ and $\mathcal P_w(w,v,z)=w$, and set
\[
 \mathcal T_S=\mathcal P_w(-\Aks)^{-1}\iota_S:\Ms\longrightarrow H_E^2(0,L).
\]

\begin{theorem}[Memory-to-beam transfer and the critical exponent $p=2$]
\label{thm:memory-beam-transfer}
Let $k_0\in H_E^2(0,L)$ be defined by $a_0(k_0,\phi)=\overline{\phi(L)}$ for all $\phi\in H_E^2(0,L)$. Then
\[
 k_0(x)=\frac{a}{2}x^2-\frac16x^3,\quad\text{with}\quad
 a=\frac{L+\frac12\delta^{(3)}L^2}{1+\delta^{(3)}L}\quad\text{and}\quad
 k_0(L)=\frac{L^3(4+\delta^{(3)}L)}{12(1+\delta^{(3)}L)},
\]
and, for every $\vartheta\in\Ms$,
\begin{equation}
 \mathcal T_S\vartheta=-\gamma\,k_0\int_0^S[\kappa(r)-\kappa(S)]\vartheta(r)\,dr.
 \label{eq:transfer-formula}
\end{equation}
Consequently
\begin{equation}
 \|\mathcal T_S\|_{\mathcal L(\Ms,H_E^2)}
 =\gamma\sqrt{k_0(L)}\,\Theta(S)^{1/2},\quad\text{with}\quad 
 \Theta(S)=\int_0^S\frac{[\kappa(r)-\kappa(S)]^2}{m(r)}\,dr,
 \label{eq:transfer-norm}
\end{equation}
where $H_E^2(0,L)$ carries the norm $\|\cdot\|_{a_0}$, and, as $S\to\infty$,
\begin{equation}
 \Theta(S)\sim
 \begin{cases}
 \dfrac{bp}{(2-p)(p+2)}\,S^{2-p}, & 1<p<2,\\[2mm]
 \dfrac b2\log S, & p=2,\\[2mm]
 \dfrac{b}{p(p-2)}, & p>2.
 \end{cases}
 \label{eq:transfer-asymptotics}
\end{equation}
Since $k_0$ is a cubic polynomial satisfying the clamped conditions, it belongs to every cubic Hermite space $V_{h_x}$, and the discrete static beam problem reproduces the profile $k_0$ exactly. This statement concerns the spatial profile. The discrete history moment and hence the scalar amplitude are generally not reproduced exactly at a fixed history mesh size.
\end{theorem}

\begin{proof}
Solve $-\Aks(w,v,z)=(0,0,\vartheta)$. Then $v=0$, $z'=\vartheta$ and $z(0)=0$, so that
\[
 J_Sz=\int_0^Sm(s)\int_0^s\vartheta(r)\,dr\,ds
 =\int_0^S[\kappa(r)-\kappa(S)]\vartheta(r)\,dr.
\]
The beam equation reduces to $a_0(w,\phi)=-\gamma J_Sz\,\overline{\phi(L)}$, hence $w=-\gamma(J_Sz)k_0$, which is \eqref{eq:transfer-formula}. The profile $k_0$ satisfies $k_0^{(4)}=0$, $k_0(0)=k_0'(0)=0$, $k_0''(L)+\delta^{(3)}k_0'(L)=0$ and $-k_0'''(L)=1$, which gives the explicit formula. With the convention that the inner product is linear in its first argument, the right-hand side of the defining identity is $\overline{\phi(L)}$. Since $\|k_0\|_{a_0}^2=a_0(k_0,k_0)=k_0(L)$ and the norm of $\vartheta\mapsto\int_0^Sg\vartheta$ on $\Ms$ is $\|g/m\|_{\Ms}$, \eqref{eq:transfer-norm} follows, with equality for $\vartheta$ proportional to $[\kappa-\kappa(S)]/m$. With $T=1+S$, direct integration gives
\[
 \Theta(S)=\frac bp\left[\frac{T^{2-p}-1}{2-p}-T^{-p}(T^2-1)+\frac{T^{2-p}-T^{-2p}}{p+2}\right],\quad\text{for }p\ne2,
\]
and $\Theta(S)=\frac b2\bigl[\log T-(1-T^{-2})+\frac14(1-T^{-4})\bigr]$ for $p=2$. Letting $S\to\infty$ yields \eqref{eq:transfer-asymptotics}. Finally, the discrete static problem is the Galerkin projection onto $V_{h_x}$ of a problem whose solution lies in $V_{h_x}$.
\end{proof}

Thus $p=2$ separates the unbounded and bounded static beam responses. The norm of the full inverse still grows as $S$ for every $p>1$.

\begin{proposition}[Convergence of the discrete static history-to-beam transfer]
\label{prop:static-transfer-convergence}
Fix $S>2$ and $\vartheta\in\mathcal M_S$. Let $z\in H^1(0,S)$ solve
\[
 z'(s)=\vartheta(s)\quad\text{and}\quad z(0)=0,
\]
and let $z_h\in V_h^s$ solve
\[
 (-\mathcal G_{S,h})z_h=P_h^m\vartheta.
\]
Then $z_h\to z$ in $\mathcal M_S$ as $h\to0$. If $z\in H^2(0,S)$, there is $C_S>0$, independent of the history mesh and of the beam mesh, such that
\begin{equation}
 \|z_h-z\|_{\mathcal M_S}\le C_S h\|z\|_{H^2(0,S)},
 \label{eq:static-history-convergence}
\end{equation}
for $0<h\le1$. For every cubic Hermite beam space $V_{h_x}$, let $w_{h_x,h}$ be the displacement component of
$(-\mathbb A_{\kappa,S}^{h_x,h})^{-1}(0,0,P_h^m\vartheta)$. Then
\[
 w_{h_x,h}=-\gamma k_0 J_Sz_h
 \quad\text{and}\quad
 \|w_{h_x,h}-\mathcal T_S\vartheta\|_{a_0}
 \le\gamma\sqrt{b k_0(L)}\,\|z_h-z\|_{\mathcal M_S}.
\]
In particular, the discrete static transfer converges for fixed $S$, with no restriction on the beam mesh. The constant $C_S$ depends on $S$.
\end{proposition}
\begin{proof}
For every $\eta_h\in V_h^s$, integration by parts and $-m'/m\ge(p+1)/(1+S)$ give
\[
 \Real(-\mathcal G_{S,h}\eta_h,\eta_h)_{m,S}
 \ge d_S\|\eta_h\|_{\mathcal M_S}^2,\quad\text{with}\quad d_S=\frac{p+1}{2(1+S)}.
\]
Thus $\|(-\mathcal G_{S,h})^{-1}\|\le d_S^{-1}$, uniformly in $h$ for fixed $S$. Assume first $z\in H^2(0,S)$ and denote its nodal interpolant by $I_hz\in V_h^s$. With $e_h=z_h-I_hz$, the defining equation and $\vartheta=z'$ imply
\[
 (-\mathcal G_{S,h})e_h
 =P_h^m\bigl(z'-(I_hz)'\bigr).
\]
Consequently,
\[
 \|z_h-z\|_{\mathcal M_S}
 \le\|I_hz-z\|_{\mathcal M_S}
 +d_S^{-1}\|(I_hz)'-z'\|_{\mathcal M_S}.
\]
Since $m$ is positive and continuous on $[0,S]$, the weighted and unweighted $L^2$ norms are equivalent for fixed $S$. The standard piecewise-linear interpolation estimates now give \eqref{eq:static-history-convergence}.

For arbitrary $z\in H^1(0,S)$ with $z(0)=0$, choose smooth $z^{(j)}$ with $z^{(j)}(0)=0$ converging to $z$ in $H^1(0,S)$. Applying the preceding estimate to each fixed $z^{(j)}$, and using the uniform discrete inverse bound $d_S^{-1}$, yields $z_h\to z$ in $\mathcal M_S$ by a density argument.

For the static coupled problem with right-hand side $(0,0,P_h^m\vartheta)$, the discrete beam velocity vanishes and its history is $z_h$. The discrete beam equation is
\[
 a_0(w_{h_x,h},\phi_{h_x})
 =-\gamma J_Sz_h\,\overline{\phi_{h_x}(L)}.
\]
The profile $k_0$ belongs to every $V_{h_x}$ and satisfies
$a_0(k_0,\phi_{h_x})=\overline{\phi_{h_x}(L)}$. Uniqueness gives $w_{h_x,h}=-\gamma(J_Sz_h)k_0$. Subtracting \eqref{eq:transfer-formula} and using $|J_Sf|\le\sqrt b\|f\|_{\mathcal M_S}$ and $\|k_0\|_{a_0}=\sqrt{k_0(L)}$ proves the estimate.
\end{proof}

\paragraph{Consistency test for the discrete static transfer.}
To test Proposition~\ref{prop:static-transfer-convergence} independently, fix $S=8$, $p=2$, $b=1$, and take $\vartheta(s)=\cos s$, so the exact transport solution is $z(s)=\sin s$. We assemble the weighted continuous $\mathbb P_1$ transport problem with twelve-point Gauss--Legendre quadrature per cell. The $L_m^2(0,S)$ errors and the errors in the moment $J_S$ are
\[
\begin{array}{c|cc}
 h & \|z_h-z\|_{\mathcal M_S} & |J_S(z_h-z)|\\\hline
 1     &2.5887\times10^{-2}&1.2797\times10^{-2}\\
 1/2   &5.4418\times10^{-3}&1.1818\times10^{-3}\\
 1/4   &1.2978\times10^{-3}&5.1085\times10^{-5}\\
 1/8   &3.2160\times10^{-4}&7.9539\times10^{-6}\\
 1/16  &8.0260\times10^{-5}&3.4303\times10^{-6}
\end{array}
\]
The history errors decrease at an empirical rate near two for this smooth datum, faster than the first-order bound of Proposition~\ref{prop:static-transfer-convergence}. The static beam profile is represented exactly by cubic Hermite elements, and the amplitude error is $\gamma\sqrt{k_0(L)}\,|J_S(z_h-z)|$.

\begin{figure}[htbp]
 \centering
 \includegraphics[width=0.94\textwidth]{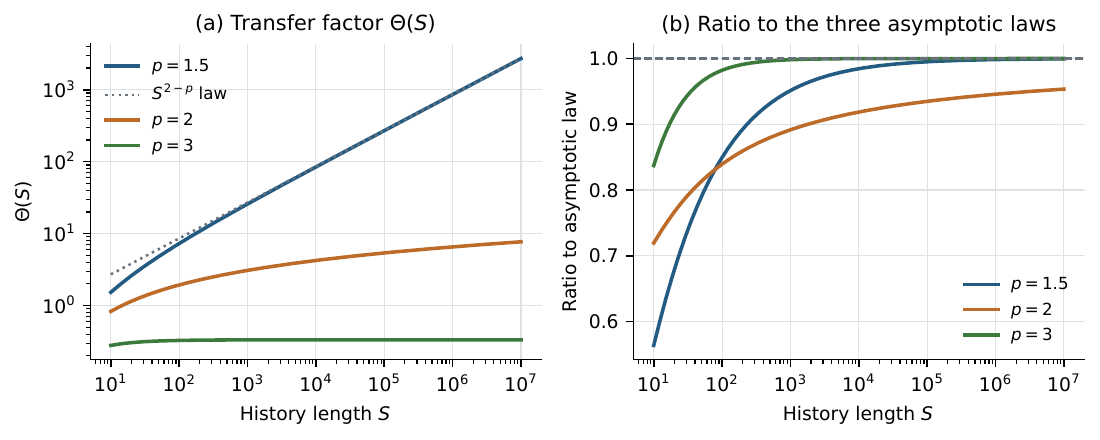}
 \caption{Memory-to-beam transfer factor for $b=1$. Left: $\Theta(S)$ for $p=1.5$, $2$, $3$ and the law $S^{2-p}$ for $p=1.5$. Right: ratio of $\Theta(S)$ to the asymptotic laws $\frac{bp}{(2-p)(p+2)}S^{2-p}$, $\frac b2\log S$ and $\frac{b}{p(p-2)}$.}
 \label{fig:memory-beam-transfer}
\end{figure}

Figure~\ref{fig:memory-beam-transfer} illustrates $\Theta(S)$ for $b=1$ and $p\in\{1.5,2,3\}$. The values follow from the explicit formulas and agree with adaptive quadrature. For $S=10^2,10^4,10^6$, the ratios to \eqref{eq:transfer-asymptotics} are $0.850$, $0.985$ and $0.998$ when $p=1.5$. They are $0.839$, $0.919$ and $0.946$ when $p=2$, and $0.982$, $1.000$ and $1.000$ when $p=3$. Convergence is slower at $p=2$ because of the additive constant in the logarithmic law.

\section{Conditional high-frequency phenomenon}
\label{app:high-frequency}
For standard semidiscretizations without distributed damping, uniform exponential bounds can fail because of insufficient high-frequency observability. Adding numerical viscosity is one way to restore it \cite{AbdallahNicaiseValeinWehbe2013}. The following proposition gives a conditional obstruction for the present scheme.

\begin{proposition}[High-frequency approximate eigenvectors of the Galerkin generator]
\label{prop:hf-quasimodes}
Let $\zeta=\delta^{(1)}=\delta^{(2)}=\delta^{(3)}=0$. Let $M_{h_x}$ and $K_{h_x}$ be the Hermite mass and stiffness matrices, $e_{h_x}$ the nodal vector with $e_{h_x}^Tq=q(L)$, and let $K_{h_x}q=\omega^2M_{h_x}q$ with $q^*M_{h_x}q=1$ and $\omega>0$. The unit vector $U=\bigl(q/(\sqrt2\,\omega),\,iq/\sqrt2,\,0\bigr)\in\Xks^{h_x,h}$ satisfies
\begin{equation*}
 \|(i\omega-\Aks^{h_x,h})U\|^2
 =\frac{|q(L)|^2}{2}\Bigl[\gamma^2e_{h_x}^*M_{h_x}^{-1}e_{h_x}+\gamma\|P_h^m1\|_{\Ms}^2\Bigr]
 \le C\,|q(L)|^2\bigl(1+h_x^{-1}\bigr)
\end{equation*}
on quasi-uniform beam meshes. Consequently, if a sequence of discrete beam eigenvectors satisfies $|q(L)|=o(h_x^{1/2})$, no bound $\|e^{t\Aks^{h_x,h}}\|\le Me^{-\omega_0t}$ holds with $M,\omega_0>0$ independent of $h_x$, even for fixed $S$ and $h$.
\end{proposition}

\begin{proof}
The first residual component vanishes. The conservative beam terms cancel because $K_{h_x}q=\omega^2M_{h_x}q$. The acceleration residual is $i\gamma q(L)M_{h_x}^{-1}e_{h_x}/\sqrt2$. The history residual is $-iq(L)P_h^m1/\sqrt2$. Their squared energy norms give the identity. The inverse trace inequality yields $e_{h_x}^*M_{h_x}^{-1}e_{h_x}\le Ch_x^{-1}$, while $\|P_h^m1\|_{\Ms}^2\le b$. Uniform exponential stability would bound $\sup_\beta\|(i\beta-\Aks^{h_x,h})^{-1}\|$ independently of $h_x$. This contradicts the vanishing residuals.
\end{proof}

The condition requires the traces and frequencies to belong to the same sequence of conservative beam eigenvectors.

\section*{Statements and declarations}

\noindent\textbf{Funding.} This research received no external funding.

\medskip
\noindent\textbf{Competing interests.} The author declares that there are no competing interests.

\end{document}